\documentclass{elsart}%
\usepackage{amsmath}
\usepackage{graphicx}
\usepackage{float}
\usepackage{indentfirst}
\usepackage{verbatim}
\usepackage{algorithm}
\usepackage{algpseudocode}
\usepackage{float}
\usepackage{multicol}
\usepackage{listings}
\usepackage{subfigure}
\usepackage{epsfig}
\usepackage{amsfonts}
\usepackage{amssymb}
\usepackage{epstopdf}
\usepackage{color}
\usepackage{epstopdf}%
\begin{document}
\begin{frontmatter}
\title{Accurate and Efficient Nystr\"{o}m Volume Integral Equation Method for the Maxwell equations for Multiple 3-D Scatterers}
\author[UNCC]{Duan Chen}
\author[UNCC]{Wei Cai}
\author[UNCC]{Brian Zinser}
\author[uml]{Min Hyung Cho}
\address[UNCC]{Department of Mathematics and Statistics, University of
North Carolina at Charlotte, Charlotte, NC 28223, USA}
\address[uml]{Department of Mathematical Sciences, University of Massachusetts Lowell, Lowell, MA 01854, USA}
\bigskip
{\bf Suggested Running Head:}
\\
Accurate and Efficient Nystr\"{o}m Volume Integral Equations for the Maxwell equations
\\
\bigskip
{\bf Corresponding Author: }
\\
Prof. Wei Cai \\
Department of Mathematics and Statistics, \\
University of North Carolina at
Charlotte, \\
Charlotte, NC 28223-0001 \\
Phone: 704-687-0628, Fax: 704-687-6415, \\
Email: wcai@uncc.edu
\newpage
\begin{abstract}
In this paper, we develop an accurate and efficient Nystr\"{o}m  volume integral equation (VIE) method for the Maxwell equations for large number of 3-D scatterers. The Cauchy Principal Values that arise from the VIE are computed accurately using a finite size exclusion volume together with explicit correction integrals consisting of removable singularities.  Also, the hyper-singular integrals are computed using interpolated quadrature formulae with tensor-product quadrature nodes for several objects, such as cubes and spheres, that are frequently encountered in the design of meta-materials . The resulting Nystr\"{o}m VIE method is shown to have high accuracy with a minimum number of collocation points and demonstrate $p$-convergence for computing the electromagnetic scattering of these objects. Numerical calculations of multiple scatterers of cubic and spherical shapes validate the efficiency and accuracy of the proposed method.
\end{abstract}
\begin{keyword}
	Electromagnetic (EM) scattering, volume integral equation, Cauchy Principal Value, dyadic Green's function, Nystr\"{o}m methods.
\end{keyword}
\textsl{AMS Subject classifications: 65R20, 65Z05, 78M25}
\end{frontmatter}

\bigskip

\section{Introduction}

Electromagnetic (EM) wave scattering of random microstructures occurs in a
wide range of applications. For example, the interaction of light with surface
plasmons on roughened metallic surfaces produces surface plasmon polaritons
(SPP) \cite{Atwater:2007,Raether:1988}, which has important applications in
solar cells \cite{Atwater:2010}, meta-materials, and super-resolution imaging
devices \cite{Pendry:2000,Fu:2010}. Also, surface enhanced Raman scattering
(SERS) \cite{Raman:1928} is closely related to the excitation of surface
plasmons on rough or nano-pattern surfaces by incident light and is a very
useful tool in finger-printing chemical components of a molecule, single
molecule detector, DNA detection, and bio-sensor, etc \cite{Hering:2008}. In
all these applications, it is critical to have accurate and efficient
numerical methods for computer simulations of the EM scattering of a large
number of microscopic objects such as spheres, cubes, etc.

In this paper, we will present an accurate and efficient Nystr\"{o}m volume
integral equation (VIE) method for the time harmonic Maxwell equations using
dyadic Green's functions $\overline{\mathbf{G}}_{\mathbf{E}}(\mathbf{r}%
,\mathbf{r}^{\prime})$. In most applications, the scatterers are embedded in
either a homogeneous or a layered medium. Thus, a VIE can be derived for the
regions occupied by the scatterers using the dyadic Green's function
$\overline{\mathbf{G}}_{\mathbf{E}}(\mathbf{r},\mathbf{r}^{\prime})$. The
resulting VIE is a second kind Fredholm integral equation and $\overline
{\mathbf{G}}_{\mathbf{E}}(\mathbf{r},\mathbf{r}^{\prime})$ will ensure that
the scattered field, expressed in terms of equivalent current sources inside
the scatterer, satisfies the Silver-M\"{u}ller radiation conditions at
infinity \cite{wcai:2013}. In the VIE formulation, computing the electric
field inside the scatterer will involve the use of Cauchy Principal Values
(CPV or simply p.v.) associated with the dyadic Green's function, namely,
\begin{equation}
{\small \text{ p.v.}\int_{\Omega}\mathrm{d}\mathbf{r}^{\prime}\text{
}\mathrm{i}\omega\overline{\mathbf{G}}_{\mathbf{E}}(\mathbf{r},\mathbf{r}%
^{\prime})\cdot\Delta\epsilon(\mathbf{r}^{\prime})\mathbf{E}(\mathbf{r}%
^{\prime})=\lim_{\delta\rightarrow0}\int_{\Omega\backslash V_{\delta}%
}\mathrm{d}\mathbf{r}^{\prime}\text{ }\mathrm{i}\omega\overline{\mathbf{G}%
}_{\mathbf{E}}(\mathbf{r},\mathbf{r}^{\prime})\cdot\Delta\epsilon
(\mathbf{r}^{\prime})\mathbf{E}(\mathbf{r}^{\prime}),} \label{cpv}%
\end{equation}
where $\Omega$ is the volume of the scatterers, $\Delta\epsilon(\mathbf{r})$
describes the scatterers' dielectric constant difference from the background
material, and $V_{\delta}$ is a small exclusion volume with size $\delta$
centered at ${\small \mathbf{r}}$. Therefore, one of the most difficult issues
for VIE methods is how to compute \textit{accurately and efficiently }the CPV
for the dyadic Green's function with an $O\left(  \displaystyle{\frac{1}%
{R^{3}}}\right)  $ singularity, where $R$ is the distance between the source
and field points.

As the CPV is defined through a limiting process of diminishing size $\delta$
of the exclusion volume $V_{\delta}{\small \mathbf{,}}$ in practical
computation, a small finite $\delta>0$ has to be taken. Then, a simple and
naive approximation could be written as
\begin{equation}
\text{ p.v.}\int_{\Omega}\mathrm{d}\mathbf{r}^{\prime}\text{ }%
{\small \mathrm{i}\omega\overline{\mathbf{G}}_{\mathbf{E}}(\mathbf{r}%
,\mathbf{r}^{\prime})\cdot\Delta\epsilon(\mathbf{r}^{\prime})\mathbf{E}%
(\mathbf{r}^{\prime})}\approx\int_{\Omega\backslash V_{\delta}}\mathrm{d}%
\mathbf{r}^{\prime}\text{ }{\small \mathrm{i}\omega\overline{\mathbf{G}%
}_{\mathbf{E}}(\mathbf{r},\mathbf{r}^{\prime})\cdot\Delta\epsilon
(\mathbf{r}^{\prime})\mathbf{E}(\mathbf{r}^{\prime})}. \label{eqn:numericalpv}%
\end{equation}

Here, we will be faced with two issues; it is the objective of this paper to
address these two issues for the Nystr\"{o}m methods and find easily
implementable solutions:

Firstly, we need to decide the size of $\delta$ to be taken and the effect of
error by using a finite $\delta$ in the calculation of the CPV. As for any
finite $\delta,$ there will be a truncation error, which we name correction
terms,
\begin{align}
\text{ p.v.}\int_{\Omega}\mathrm{d}\mathbf{r}^{\prime}\text{ }%
{\small \mathrm{i}\omega\overline{\mathbf{G}}_{\mathbf{E}}(\mathbf{r}%
,\mathbf{r}^{\prime})\cdot\Delta\epsilon(\mathbf{r}^{\prime})\mathbf{E}%
(\mathbf{r}^{\prime})} &  =\int_{\Omega\backslash V_{\delta}}\mathrm{d}%
\mathbf{r}^{\prime}\text{ }{\small \mathrm{i}\omega\overline{\mathbf{G}%
}_{\mathbf{E}}(\mathbf{r},\mathbf{r}^{\prime})\cdot\Delta\epsilon
(\mathbf{r}^{\prime})\mathbf{E}(\mathbf{r}^{\prime})}\nonumber\\
&  +\text{ correction terms}.\label{correction}%
\end{align}

Consequently, the correction terms will play an important role in determining
the accurate numerical solution. The correction terms were derived by Fikioris
\cite{Fik:1965} using a mixed potential formulation of the electric field. In
this paper, we will re-derive the VIE using vector and scalar potentials such
that the CPV can be computed with explicit expressions for the correction
terms and the resulting VIE is in a form more suitable for numerical
implementations. Previous work on how to handle singular integrals for VIE
methods include singularity subtraction \cite{Kottmann:2000}, locally
corrected Nystr\"{o}m scheme \cite{Liu:2001}, direct integration of the
singularity \cite{Tong:2010},\ etc.

Secondly, we need to find accurate and efficient ways to numerically compute
the integral over $\Omega\backslash V_{\delta}$ for a fixed exclusion volume
$V_{\delta}$.

To address the second issue, we will use special quadrature weights over
tensor-product quadrature nodes in a reference element $\Omega$ (a sphere or a
cube in this paper) by using an interpolation approach. Specifically, first a
brute-force computation of the integral, using Gauss quadrature in polar
coordinate centered at the singularity, will be done to a given accuracy with
a large number of evaluations of the integrand. The integrands in the VIE
matrix entries, except for the singular denominators involving $R^{k},1\leq
k\leq3$, are smooth functions, therefore they can be accurately interpolated
using values only on a small number of the tensor-product nodes inside the
domain $\Omega.$ Then , the brute-force integration formula weights can be
converted into new integration weights for the tensor-product nodes. Moreover,
the new integration formula can be tabulated for integrating general
functions. The Nystr\"{o}m collocation method based on the simple
tensor-product nodes will be used to solve the VIE for the scattering of a
large number of scatterers with a high accuracy using a minimal number of unknowns.

The rest of the paper is organized as follows: Section 2 presents the
formulation of a VIE where the CPV can be computed with a finite exclusion
volume $V_{\delta}$ accurately. Then, numerical algorithms are given in
Section 3: the Nystr\"{o}m collocation methods, an efficient quadrature
formula, and numerical implementation. Section \ref{sec:results} includes
various numerical tests such as the accuracy of CPV computation, $\delta
$-independence of matrix entries, $p$-refinement convergence for the VIE for
spheres and cubes, and results of scattering from multiple scatterers. The
paper ends with a conclusion in Section \ref{sec:conclusion}.

\section{Volume Integral equations for Maxwell equations}

%\bigskip

\subsection{VIE and Cauchy principal values}

In this section, we will follow \cite{wcai:2013} to show briefly how the VIE
for the Maxwell equations can be derived using a vector form of Green's second
identity for the following vector wave equations for the electric field
$\mathbf{E}(\mathbf{r})$,
\begin{equation}
\mathcal{L}\mathbf{E}(\mathbf{r})\text{ }\mathbf{-}\text{ }\omega^{2}%
\epsilon(\mathbf{r})\mathbf{E}(\mathbf{r})=-\mathrm{i}\omega\mathbf{J}%
_{e}(\mathbf{r}),\text{ \ \ }\mathbf{r}\in\mathbb{R}^{3}\backslash
\partial\Omega, \label{per-171}%
\end{equation}
where $\omega$ is the frequency, $\mu$ is the permeability, and $\mathcal{L}$
is a differential operator defined by
\[
\mathcal{L}=\nabla\times\frac{1}{\mu}\nabla\times,
\]
source $\mathbf{J}_{e}(\mathbf{r})$ produces the incident wave impinging on
the scatterers from above, i.e.,%
\begin{equation}
\mathbf{E}^{\text{inc}}(\mathbf{r})=-\mathrm{i}\omega\mu(\mathbf{r}%
)\int_{\mathbb{R}^{3}}\overline{\mathbf{G}}_{E}(\mathbf{r},\mathbf{r}^{\prime
})\cdot\mathbf{J}_{e}(\mathbf{r}^{\prime})\mathrm{d}\mathbf{r}^{\prime},
\label{per-173}%
\end{equation}
and $\overline{\mathbf{G}}_{E}(\mathbf{r},\mathbf{r}^{\prime})$ is the dyadic
Green's function to be defined later. A scatterer $\Omega$ is characterized by
a different dielectric constant from that of the background dielectrics
$\epsilon_{L}(\mathbf{r})$, i.e.,%
\begin{equation}
\epsilon(\mathbf{r})=\epsilon_{L}(\mathbf{r})+\Delta\epsilon(\mathbf{r}),
\label{per-175}%
\end{equation}
where $\Delta\epsilon(\mathbf{r})=0,\mathbf{r\notin}$ $\Omega.$ Then,
(\ref{per-171}) can be rewritten as%
\begin{equation}
\mathcal{L}\mathbf{E}(\mathbf{r})\text{ }\mathbf{-}\text{ }\omega^{2}%
\epsilon_{L}(\mathbf{r})\mathbf{E}(\mathbf{r})=-\mathrm{i}\omega
\mathbf{J}(\mathbf{r}), \label{per-177}%
\end{equation}
where%
\begin{equation}
\mathbf{J}(\mathbf{r})=\mathbf{J}_{e}(\mathbf{r})+\mathbf{J}_{\text{eq}%
}(\mathbf{r}), \label{per-178}%
\end{equation}
and the equivalent current source $\mathbf{J}_{\text{eq}}(\mathbf{r})$ is
defined to characterize the presence of the scatterer $\Omega$:%
\begin{equation}
\mathbf{J}_{\text{eq}}(\mathbf{r})=\mathrm{i}\omega\Delta\epsilon
(\mathbf{r})\mathbf{E}(\mathbf{r}). \label{per-179}%
\end{equation}

Let us consider any point $\mathbf{r}^{\prime}\in\Omega$ and a small volume
$V_{\delta}=V_{\delta}(\mathbf{r}^{\prime})\subset\Omega$ centered at
$\mathbf{r}^{\prime}.$ The dyadic Green's function $\overline{\mathbf{G}}%
_{E}(\mathbf{r},\mathbf{r}^{\prime})$ is defined by
\begin{equation}
\mathcal{L}\overline{\mathbf{G}}_{E}(\mathbf{r},\mathbf{r}^{\prime})\text{{}%
}\mathbf{-}\text{ }\omega^{2}\epsilon_{L}(\mathbf{r})\overline{\mathbf{G}}%
_{E}(\mathbf{r},\mathbf{r}^{\prime})=\frac{1}{\mu(\mathbf{r})}\overline
{\mathbf{I}}\delta(\mathbf{r}-\mathbf{r}^{\prime}),\text{ }\mathbf{r}%
\in\mathbb{R}^{3}. \label{per-185}%
\end{equation}
In case of free-space, we have%
\begin{equation}
\overline{\mathbf{G}}_{\mathbf{E}}\mathbf{(r},\mathbf{r^{\prime}%
)=\overline{\mathbf{G}}_{E}(r^{\prime}},\mathbf{r)=}\left(  \overline
{\mathbf{I}}\text{ }\mathbf{+}\text{ }\frac{1}{k^{2}}\nabla\nabla\right)
g(\mathbf{r},\mathbf{r^{\prime}}), \label{eqn3-34}%
\end{equation}
where $k^{2}=\omega^{2}\epsilon_{L}(\mathbf{r})\mu$ and
\begin{equation}
g(\mathbf{r},\mathbf{r^{\prime}})=\frac{1}{4\pi}\frac{\mathrm{e}%
^{-\mathrm{i}kR}}{R},\qquad R=|\mathbf{r}-\mathbf{r}^{\prime}|.
\label{eqn3-35}%
\end{equation}

Next, we multiply (\ref{per-177}) by $\overline{\mathbf{G}}_{E}(\mathbf{r}%
,\mathbf{r}^{\prime})$ and (\ref{per-185}) by $\mathbf{E}(\mathbf{r})$, form
the difference, and integrate over the domain $\mathbb{R}^{3}\backslash
V_{\delta}.$ After some manipulation \cite{wcai:2013}, we arrive at the
following equation (after switching $\mathbf{r}$ and $\mathbf{r}^{\prime}$):%
\begin{align}
&  -\mathrm{i}\omega\mu(\mathbf{r})\int_{\mathbb{R}^{3}\backslash V_{\delta}%
}\mathrm{d}\mathbf{r}^{\prime}\text{ }\overline{\mathbf{G}}_{E}(\mathbf{r}%
,\mathbf{r}^{\prime})\cdot\mathbf{J}(\mathbf{r}^{\prime})-\mu(\mathbf{r}%
)\int_{{S}_{\delta}}\ \mathrm{d}s^{\prime}\left[  \mathrm{i}\omega\text{{}%
}\overline{\mathbf{G}}_{E}(\mathbf{r},\mathbf{r}^{\prime})\cdot\left(
\mathbf{n}\times\mathbf{H}(\mathbf{r}^{\prime})\right)  \right. \nonumber\\
&  {-}\text{ }\frac{1}{\mu(\mathbf{r}^{\prime})}\nabla\times\overline
{\mathbf{G}}_{E}(\mathbf{r},\mathbf{r}^{\prime})\cdot\left(  \mathbf{n\times
E}(\mathbf{r}^{\prime})\right)  ]=\mathbf{0},\qquad\qquad\mathbf{r}\in\Omega,
\label{per-186}%
\end{align}
where $S_{\delta}=\partial V_{\delta}(\mathbf{r}),$ $\mathbf{n}$ is the normal
of $S_{\delta}$ pointing out of $V_{\delta}(\mathbf{r})$.

As $\delta\rightarrow0$, the first integral will approach the CPV of a
singular integral, while the surface integrals depends on the geometric shape
of the volume $V_{\delta}.$

In order to estimate the surface integrals, we have the following asymptotic
approximations for $kR\ll1$\cite{Bladel}:
\begin{align}
\overline{\mathbf{G}}_{\mathbf{E}}\mathbf{(r},\mathbf{r^{\prime})}  &
=\frac{1}{4\pi k^{2}R^{3}}(\mathbf{I}-3\mathbf{u\otimes u})+O\left(  \frac
{1}{R^{2}}\right)  ,\label{per-187}\\
\nabla^{\prime}\times\overline{\mathbf{G}}_{\mathbf{E}}\mathbf{(r}%
,\mathbf{r^{\prime})}  &  =\frac{1}{4\pi R^{2}}\mathbf{u}\times\mathbf{I}%
+O\left(  \frac{1}{R}\right)  , \label{per-189}%
\end{align}
where $\mathbf{u}={(\mathbf{r}}^{\prime}{-\mathbf{r})}/{R},$ which implies
that:
\begin{align}
&  \underset{\delta\rightarrow0}{\text{lim}}\int_{{S}_{\delta}}\ \mathrm{d}%
s^{\prime}\text{ }\mathbf{n\times E}(\mathbf{r}^{\prime})\cdot\nabla
\times\overline{\mathbf{G}}_{E}(\mathbf{r}^{\prime},\mathbf{r})=-\left[
\mathbf{I-L}_{V_{\delta}}\right]  \cdot\mathbf{E}(\mathbf{r}), \label{per-191}%
\\
&  \underset{\delta\rightarrow0}{\text{lim}}\int_{{S}_{\delta}}\ \mathrm{d}%
s^{\prime}\text{ }\mathbf{n}\times\mathbf{H}(\mathbf{r}^{\prime}%
)\cdot\overline{\mathbf{G}}_{E}(\mathbf{r}^{\prime},\mathbf{r})=-\frac
{1}{k^{2}}\mathbf{L}_{V_{\delta}}\cdot\nabla\times\mathbf{H}(\mathbf{r}),
\label{per-193a}%
\end{align}
and the $\mathbf{L}$-dyadics for\ $V_{\delta}$ of various geometric shapes are
given in \cite{Yaghjian:1980}. We have
\begin{equation}
\mathbf{L}_{V_{\delta}}=\frac{1}{3}\mathbf{I}%
\end{equation}
for a sphere as used in this paper.

Substituting (\ref{per-191}) and (\ref{per-193a}) into (\ref{per-186}), and
using Amp\`{e}re's law, we have the VIE for the electric field for
$\mathbf{r}\in\Omega$:
\begin{equation}
\mathbf{C\cdot E}(\mathbf{r})=\mathbf{E}^{\text{inc}}(\mathbf{r}%
)-\mathrm{i}\omega\mu(\mathbf{r})\text{ p.v.}\int_{\Omega}\mathrm{d}%
\mathbf{r}^{\prime}\text{ }{\small \mathrm{i}\omega\overline{\mathbf{G}%
}_{\mathbf{E}}(\mathbf{r},\mathbf{r}^{\prime})\cdot\Delta\epsilon
(\mathbf{r}^{\prime})\mathbf{E}(\mathbf{r}^{\prime})}, \label{per-197}%
\end{equation}
where the coefficient matrix is given by
\begin{equation}
\mathbf{C=I+L}_{V_{\delta}}\cdot\Delta\epsilon(\mathbf{r}). \label{per-199}%
\end{equation}

As mentioned in the introduction, in most numerical implementations of
(\ref{per-197}), the CPV integral is computed by selecting a finite $\delta$
as in (\ref{eqn:numericalpv}). The numerical solution thus obtained does not
satisfy the original VIE (\ref{per-197}) and will have an intrinsic error
related to the correction terms in (\ref{correction}).

\subsection{Reformulation of the VIE and computing CPV with a finite $\delta$}

\label{sec:newvie}

%In the numerical computation of (), a the p.v. is generally approximated by
%selecting a finite and small $\delta,$ namely,%
%\begin{equation}
%\mathbf{C\cdot E}(\mathbf{r})=\mathbf{E}^{\text{inc}}(\mathbf{r}%
%)-\mathrm{i}\omega\mu(\mathbf{r})\text{ }\int_{\Omega\backslash V_{\delta}%
%}\mathrm{d}\mathbf{r}^{\prime}\text{ }\mathrm{i}\omega\Delta\epsilon
%(\mathbf{r}^{\prime})\mathbf{E}(\mathbf{r}^{\prime})\cdot\overline{\mathbf{G}%
%}_{E}(\mathbf{r}^{\prime},\mathbf{r})+O(\delta). \label{vie001}%
%\end{equation}
%Therefore, the numerical solution is in fact an approximation to () by
%ignoring the $O(\delta)$ - which is a "truncation error" of the IE, implying
%that the accuracy of the numerical solution will be always limited by this
%truncation error on the level of IE no matter how accurately the other terms
%of the IE (\ref{vie001}) is evaluated. In this section, we will re-derive the
%VIE (\ref{per-197}) by using vector-scalar potentials for the electric field
%and give the explicit formula for the truncation error in (\ref{vie001}).

In this section, we will re-derive the VIE for the electric field where the
CPV in (\ref{per-199}) can be computed with a finite exclusion volume
$V_{\delta}$ together with some correction terms. Based on the Helmholtz
decomposition, the electric field $\mathbf{E(r)}$ can be expressed as follows
\cite{wcai:2013},%
\[
\mathbf{E}=-\mathrm{i}\omega\mathbf{A}-\nabla V,
\]
where $A$ and $V$ are vector and scalar potentials, respectively. Using the
Lorentz gauge condition \cite{stratton:1941},
\begin{equation}
\nabla\cdot\mathbf{A}=-\mathrm{i}\omega\epsilon\mu V, \label{eqn1-45}%
\end{equation}
we can have a vector potential representation
%\begin{equation}
%\mathbf{J}(\mathbf{r})=\mathbf{J}_{e}(\mathbf{r})+\mathbf{J}_{\text{eq}%
%}(\mathbf{r}), \label{per-178}%
%\end{equation}%
\begin{equation}
\mathbf{E}=-\mathrm{i}\omega\mathbf{A}+\frac{1}{\mathrm{i}\omega\epsilon\mu
}\nabla(\nabla\cdot\mathbf{A})=-\mathrm{i}\omega\left[  \overline{\mathbf{I}%
}\text{ }\mathbf{+}\text{ }\frac{1}{k^{2}}\nabla\nabla\right]  \mathbf{A}.
\label{eqn1-49}%
\end{equation}

Meanwhile, it can be shown that the potential $\mathbf{A}$ satisfies the
Helmholtz equation component-wise \cite{stratton:1941},
\begin{equation}
\nabla^{2}\mathbf{A}+k^{2}\mathbf{A}=-\mu\mathbf{J}. \label{eqn:1-47}%
\end{equation}

Thus, the solution $\mathbf{A}$ of Eq. (\ref{eqn:1-47}) can be rewritten in an
integral representation:
\begin{align}
\mathbf{A}  &  =\mu\int_{\mathbb{R}^{3}}\mathrm{d}\mathbf{r}^{\prime
}\mathbf{J}(\mathbf{r}^{\prime})g(\mathbf{r},\mathbf{r}^{\prime}%
)\nonumber\label{per-180}\\
&  =\mu\int_{\mathbb{R}^{3}\backslash\Omega}\mathbf{J}_{e}(\mathbf{r}^{\prime
})g(\mathbf{r},\mathbf{r}^{\prime})\mathrm{d}\mathbf{r}^{\prime}+\mu
\int_{\Omega}\mathbf{J}_{\mathrm{eq}}(\mathbf{r}^{\prime})g(\mathbf{r}%
,\mathbf{r}^{\prime})\mathrm{d}\mathbf{r}^{\prime}\nonumber\\
&  =\mu\int_{\mathbb{R}^{3}\backslash\Omega}\mathbf{J}_{e}(\mathbf{r}^{\prime
})g(\mathbf{r},\mathbf{r}^{\prime})\mathrm{d}\mathbf{r}^{\prime}+\mu
\int_{\Omega}\mathrm{d}\mathbf{r}^{\prime}\mathrm{i}\omega\Delta
\epsilon(\mathbf{r}^{\prime})\mathbf{E}(\mathbf{r}^{\prime})g(\mathbf{r}%
,\mathbf{r}^{\prime}),
\end{align}
where the second equality on the right hand side of Eq. (\ref{per-180}) is due
to the assumption that $\mathrm{supp}(\mathbf{J}_{e}(\mathbf{r}))\cap
\Omega=\emptyset.$

The first integral of (\ref{per-180}) is well-defined if $\mathbf{r}\in\Omega$
and after being plugged it into Eq. (\ref{eqn1-49}), it yields the incident
wave $\mathbf{E}^{\mathrm{inc}}(\mathbf{r})$ according to relation
(\ref{per-173}). For the second integral over $\Omega$ , we split it as
follows:
\begin{equation}
\mu\int_{\Omega}\mathrm{d}\mathbf{r}^{\prime}\mathrm{i}\omega\Delta
\epsilon(\mathbf{r}^{\prime})\mathbf{E}(\mathbf{r}^{\prime})g(\mathbf{r}%
,\mathbf{r}^{\prime})=\mu\left(  \int_{\Omega\backslash V_{\delta}}%
+\int_{V_{\delta}}\right)  \mathrm{i}\omega\Delta\epsilon(\mathbf{r}^{\prime
})\mathbf{E}(\mathbf{r}^{\prime})g(\mathbf{r},\mathbf{r}^{\prime}),\nonumber
\end{equation}
and along with the first integral term in Eq. (\ref{per-180}), it follows from
Eq. (\ref{eqn1-49}) that%
\begin{align}
\mathbf{E}  &  =\mathbf{E}^{\mathrm{inc}}(\mathbf{r})-\mathrm{i}\omega\mu
\int_{\Omega\backslash V_{\delta}}\mathrm{i}\omega\Delta\epsilon
(\mathbf{r}^{\prime})\mathbf{E}(\mathbf{r}^{\prime})\left[  \overline
{\mathbf{I}}\mathbf{+}\frac{1}{k^{2}}\nabla\nabla\right]  g(\mathbf{r}%
,\mathbf{r}^{\prime})\nonumber\label{per-207}\\
&  -\mathrm{i}\omega\mu\left[  \overline{\mathbf{I}}\mathbf{+}\frac{1}{k^{2}%
}\nabla\nabla\right]  \int_{V_{\delta}}\mathrm{d}\mathbf{r}^{\prime}%
\mathrm{i}\omega\Delta\epsilon(\mathbf{r}^{\prime})\mathbf{E}(\mathbf{r}%
^{\prime})g(\mathbf{r},\mathbf{r}^{\prime})\nonumber\\
&  =\mathbf{E}^{\mathrm{inc}}(\mathbf{r})-\mathrm{i}\omega\mu\int
_{\Omega\backslash V_{\delta}}\mathrm{i}\omega\Delta\epsilon(\mathbf{r}%
^{\prime})\overline{\mathbf{G}}_{\mathbf{E}}\mathbf{(r},\mathbf{r^{\prime
})\cdot E}(\mathbf{r}^{\prime})\nonumber\\
&  -\mathrm{i}\omega\mu\left[  \overline{\mathbf{I}}\mathbf{+}\frac{1}{k^{2}%
}\nabla\nabla\right]  \int_{V_{\delta}}\mathrm{d}\mathbf{r}^{\prime}%
\mathrm{i}\omega\Delta\epsilon(\mathbf{r}^{\prime})\mathbf{E}(\mathbf{r}%
^{\prime})g(\mathbf{r},\mathbf{r}^{\prime}).
\end{align}
Next, we separate the singular part in $g(\mathbf{r},\mathbf{r}^{\prime})$ as
follows:
\begin{equation}
g(\mathbf{r},\mathbf{r^{\prime}})=g_{0}(\mathbf{r},\mathbf{r^{\prime}%
})+\widetilde{g}(\mathbf{r},\mathbf{r^{\prime}}), \label{sigsplit}%
\end{equation}
where%
\begin{equation}
g_{0}(\mathbf{r},\mathbf{r^{\prime}})=\frac{1}{4\pi|\mathbf{r}%
-\mathbf{r^{\prime}|}},\text{ \ }\widetilde{g}=g-g_{0}. \label{per-211}%
\end{equation}
Then, using the fact that \cite{Lee:1980,Yaghjian:1980}
\begin{equation}
\nabla\nabla\int_{V_{\delta}}\mathrm{d}\mathbf{r}^{\prime}\frac{1}%
{4\pi|\mathbf{r}-\mathbf{r^{\prime}|}}=-\int_{\partial V_{\delta}}%
\mathrm{d}s^{\prime}\frac{\left(  \mathbf{r}-\mathbf{r^{\prime}}\right)
\mathbf{u}_{n}(\mathbf{r^{\prime}})}{4\pi|\mathbf{r}-\mathbf{r^{\prime}|}^{3}%
}=-\mathbf{L}_{V_{\delta}}, \label{per-213}%
\end{equation}
we can compute the following integral as%
\begin{align}
&  \nabla\nabla\int_{V_{\delta}}\mathrm{d}\mathbf{r}^{\prime}\Delta
\epsilon(\mathbf{r}^{\prime})\mathbf{E}(\mathbf{r}^{\prime})g_{0}%
(\mathbf{r},\mathbf{r}^{\prime})\nonumber\\
&  =\nabla\nabla\int_{V_{\delta}}\mathrm{d}\mathbf{r}^{\prime}\frac{1}%
{4\pi|\mathbf{r}-\mathbf{r^{\prime}|}}\Delta\epsilon(\mathbf{r})\mathbf{E}%
(\mathbf{r})+\int_{V_{\delta}}\mathrm{d}\mathbf{r}^{\prime}\nabla\nabla
g_{0}(\mathbf{r},\mathbf{r}^{\prime})\left[  \Delta\epsilon(\mathbf{r}%
^{\prime})\mathbf{E}(\mathbf{r}^{\prime})-\Delta\epsilon(\mathbf{r}%
)\mathbf{E}(\mathbf{r})\right] \nonumber\\
&  =-\mathbf{L}_{V_{\delta}}\Delta\epsilon(\mathbf{r})\mathbf{E}%
(\mathbf{r})+\int_{V_{\delta}}\mathrm{d}\mathbf{r}^{\prime}\nabla\nabla
g_{0}(\mathbf{r},\mathbf{r}^{\prime})\left[  \Delta\epsilon(\mathbf{r}%
^{\prime})\mathbf{E}(\mathbf{r}^{\prime})-\Delta\epsilon(\mathbf{r}%
)\mathbf{E}(\mathbf{r})\right]  , \label{per-215}%
\end{align}
where the second integral has a removable singularity $O\left(  {\frac
{1}{|\mathbf{r}-\mathbf{r^{\prime}|}^{2}}}\right)  $ through the use of
spherical coordinates centered at $\mathbf{r}$, provided the function
$\Delta\epsilon(\mathbf{r})\mathbf{E}(\mathbf{r})$ is differentiable in the
interior of $V_{\delta}$, which we assume it to be.

With (\ref{per-215}) and the fact that $\tilde{g}=g-g_{0}$ is a smooth
function, (\ref{per-207}) becomes%
\begin{align}
\mathbf{C\cdot E}  &  =\mathbf{E}^{\text{inc}}(\mathbf{r})-i\omega\mu
\int_{\Omega\backslash V_{\delta}}i\omega\Delta\epsilon(\mathbf{r}^{\prime
})\overline{\mathbf{G}}_{\mathbf{E}}\mathbf{(r},\mathbf{r^{\prime})\cdot
E}(\mathbf{r}^{\prime})\nonumber\\
&  +\omega^{2}\mu\mathrm{\ }\int_{V_{\delta}}\mathrm{d}\mathbf{r}^{\prime
}\Delta\epsilon(\mathbf{r}^{\prime})\mathbf{E}(\mathbf{r}^{\prime
})g(\mathbf{r},\mathbf{r}^{\prime})\nonumber\\
&  +\frac{\omega^{2}}{k^{2}}\mu\mathrm{\ }\int_{V_{\delta}}\mathrm{d}%
\mathbf{r}^{\prime}\Delta\epsilon(\mathbf{r}^{\prime})\nabla\nabla
\widetilde{g}(\mathbf{r},\mathbf{r}^{\prime})\cdot\mathbf{E}(\mathbf{r}%
^{\prime})\nonumber\\
&  +\frac{\omega^{2}}{k^{2}}\mu\int_{V_{\delta}}\mathrm{d}\mathbf{r}^{\prime
}\nabla\nabla g_{0}(\mathbf{r},\mathbf{r}^{\prime})\left[  \Delta
\epsilon(\mathbf{r}^{\prime})\mathbf{E}(\mathbf{r}^{\prime})-\Delta
\epsilon(\mathbf{r})\mathbf{E}(\mathbf{r})\right]  , \label{eqn:per-217}%
\end{align}
with same coefficient $\mathbf{C}$ as in (\ref{per-199}).

The VIE in (\ref{eqn:per-217}) is similar to those obtained by Fikioris
\cite{Fik:1965}; however, our derivation is based on a splitting of Green's
function in (\ref{sigsplit}) and the identity for $\mathbf{L}_{V_{\delta}}$ in
(\ref{per-213}). A comparison study between CPV formulation (\ref{per-197})
and finite exclusion volume formulation (\ref{eqn:per-217}) can be found in
\cite{Wang:1982}. Now expression (\ref{eqn:per-217}) holds for any finite
$\delta>0$ as long as $V_{\delta}\subset\Omega$, and all integrals involved on
the right-hand side are well-defined provided that $\Delta\epsilon
(\mathbf{r})\mathbf{E}(\mathbf{r})$ is H\"{o}lder continuous. We can see that
the last three integrals can be understood as the correction terms for
computing the Cauchy principal value with a finite-sized exclusion volume
$V_{\delta}$. It should be noted that these integrals are all weakly singular
integrals whose singularities can be removed by a spherical coordinate
transform. In particular, we can estimate their magnitudes in terms of
$\delta.$ Namely,%
\begin{equation}
|\int_{V_{\delta}}\mathrm{d}\mathbf{r}^{\prime}\Delta\epsilon(\mathbf{r}%
^{\prime})\mathbf{E}(\mathbf{r}^{\prime})g(\mathbf{r},\mathbf{r}^{\prime
})|\leq C_{1}||\Delta\epsilon\mathbf{E||}_{\infty}\delta^{2}, \label{per-219}%
\end{equation}%
\begin{equation}
|\int_{V_{\delta}}\mathrm{d}\mathbf{r}^{\prime}\Delta\epsilon(\mathbf{r}%
^{\prime})\nabla\nabla\widetilde{g}(\mathbf{r},\mathbf{r}^{\prime}%
)\cdot\mathbf{E}(\mathbf{r}^{\prime})|\leq C_{2}||\Delta\epsilon
\mathbf{E||}_{\infty}\delta^{2}, \label{per-221}%
\end{equation}
and
\begin{equation}
|\int_{V_{\delta}}\mathrm{d}\mathbf{r}^{\prime}\nabla\nabla g_{0}%
(\mathbf{r},\mathbf{r}^{\prime})\left[  \Delta\epsilon(\mathbf{r}^{\prime
})\mathbf{E}(\mathbf{r}^{\prime})-\Delta\epsilon(\mathbf{r})\mathbf{E}%
(\mathbf{r})\right]  |\leq C_{3}||\Delta\epsilon\mathbf{E||}_{1,\infty}\delta,
\label{per-223}%
\end{equation}
where $C_{1}$, $C_{2}$, $C_{3}$ are constants, and $\Vert\cdot\Vert_{\infty}$
and $\Vert\cdot\Vert_{1,\infty}$ represent the $L^{\infty}$ norms of a
function and its first derivative, respectively.

\begin{rem}
\noindent\ Equations (\ref{per-219})-(\ref{per-223}) explicitly show the
accuracy of approximating the CPV (\ref{cpv}) by the integral
(\ref{eqn:numericalpv}) with a finite $\delta>0$, i.e., the truncation error
is of the order $O(\delta)$. Hence the numerical solution of the VIE will have
this $O(\delta)$ truncation error in general regardless of the integration
quadratures used if terms in (\ref{per-219})-(\ref{per-223}) are not included.

Moreover, in case that $V_{\delta}$ is a ball of radius $\delta$ centered at
$\mathbf{r}$, one can obtain a better estimate than (\ref{per-223}) due to the
anti-symmetry of the singular term $\nabla\nabla g_{0}(\mathbf{r}%
,\mathbf{r}^{\prime})$ in the spherical coordinates, i.e.
\begin{equation}
|\int_{V_{\delta}}\mathrm{d}\mathbf{r}^{\prime}\nabla\nabla g_{0}%
(\mathbf{r},\mathbf{r}^{\prime})\left[  \Delta\epsilon(\mathbf{r}^{\prime
})\mathbf{E}(\mathbf{r}^{\prime})-\Delta\epsilon(\mathbf{r})\mathbf{E}%
(\mathbf{r})\right]  |\leq C_{4}||\Delta\epsilon\mathbf{E||}_{2,\infty}%
\delta^{2}, \label{eqn:est}%
\end{equation}
where $C_{4}$ are constant.
\end{rem}

\section{Numerical Methods}

\label{sec:methods}

\subsection{Nystr\"{o}m collocation method}

We use Nystr\"{o}m collocation methods to solve (\ref{eqn:per-217}). First, we
assume the computational domain $\Omega$ comprised of an $N-$number
non-overlapping elements (cubes and balls) $\Omega_{i}$ $,i=1,2,...,N$. On
each element $\Omega_{i}$, we assign $M$ tensor-product nodes for which $M$
scalar interpolant basis functions $\phi_{ij},j=1,2,3,....M$ are defined.
Then, we can write the solution as
\begin{equation}
\mathbf{E}({}\mathbf{r})=\sum_{i=1}^{N}\sum_{j=1}^{M}\mathbf{c}_{ij}\phi
_{ij}(\mathbf{r}),\text{ \ \ \ \ }\mathbf{r\in}\Omega_{i},
\label{eqn:discretization}%
\end{equation}
where $\mathbf{c}_{ij}$ are the unknown vectorial coefficients. Inserting
(\ref{eqn:discretization}) into (\ref{eqn:per-217}),
%\begin{equation}\nonumber
%	\phi_{ij}({\bf r}_{\tilde{i}\tilde{j}})=1 \quad\text{ only when } i=\tilde{i} \text{ and } j=\tilde{j}
%\end{equation}
we obtain the following equations for $\mathbf{c}_{ij}$:
%{\small\begin{eqnarray}\nonumber
%	\mathbf{C} \cdot \sum_{i=1}^N\sum_{j=1}^M\mathbf{c}_{ij}\phi_{ij}({\bf r}_{ij})&=&\mathbf{E}^{\rm inc}({\bf r}_{ij})+\omega^2\mu\sum_{i=1}^N\sum_{j=1}^M\left[\int_{\Omega\backslash V_{\delta_{ij}}}\mathrm{d}{\bf r}'\Delta\epsilon({\bf r}')\bar{\bf G}_{\bf E}({\bf r}_{ij}, {\bf r}')\phi_{ij}({\bf r}')\right]{\bf c}_{ij}\\\nonumber
%	&+&\omega^2\mu\sum_{i=1}^N\sum_{j=1}^M\left[\int_{V_{\delta_{ij}}}\mathrm{d}{\bf r}'\Delta\epsilon(\mathbf{r}^{\prime})g(\mathbf{r}_{ij},\mathbf{r}^{\prime})\phi_{ij}({\bf r}')\right]{\bf c}_{ij}\\\nonumber
%	&+&\frac{\omega^2\mu}{k^2_L}\sum_{i=1}^N\sum_{j=1}^M\left[\int_{V_{\delta_{ij}}}\mathrm{d}{\bf r}'\Delta\epsilon(\mathbf{r}^{\prime})\nabla\nabla
%\widetilde{g}(\mathbf{r}_{ij},\mathbf{r}^{\prime})\phi_{ij}({\bf r}')\right]\cdot{\bf c}_{ij}\\
%	&+&\frac{\omega^{2}\mu}{k_{L}^{2}}{\int_{V_{\delta_{ij}}}\mathrm{d}\mathbf{r}%
%^{\prime}\nabla\nabla g_{0}(\mathbf{r}_{ij},\mathbf{r}^{\prime})\left[
%\Delta\epsilon(\mathbf{r}^{\prime})\phi_{ij}(\mathbf{r}^{\prime}%
%)-\Delta\epsilon(\mathbf{r}_{ij})\phi_{ij}(\mathbf{r}_{ij})\right]}\cdot{\bf c}_{ij}
%\end{eqnarray}}
{\small
\begin{align}
\mathbf{C}\cdot\mathbf{c}_{ij}  &  =\mathbf{E}_{ij}^{\mathrm{inc}}+\omega
^{2}\mu\sum_{n=1}^{N}\sum_{m=1}^{M}\left[  \int_{\Omega_{n}\backslash
V_{\delta_{ij}}}\mathrm{d}\mathbf{r}^{\prime}\Delta\epsilon(\mathbf{r}%
^{\prime})\bar{\mathbf{G}}_{\mathbf{E}}(\mathbf{r}_{ij},\mathbf{r}^{\prime
})\phi_{nm}(\mathbf{r}^{\prime})\right]  \cdot\mathbf{c}_{nm}%
\nonumber\label{eqn:discre}\\
&  +\omega^{2}\mu\sum_{m=1}^{M}\left[  \int_{V_{\delta_{ij}}}\mathrm{d}%
\mathbf{r}^{\prime}\Delta\epsilon(\mathbf{r}^{\prime})g(\mathbf{r}%
_{ij},\mathbf{r}^{\prime})\phi_{im}(\mathbf{r}^{\prime})\right]
\cdot\mathbf{c}_{im}\nonumber\\
&  +\frac{\omega^{2}\mu}{k^{2}}\sum_{m=1}^{M}\left[  \int_{V_{\delta_{ij}}%
}\mathrm{d}\mathbf{r}^{\prime}\Delta\epsilon(\mathbf{r}^{\prime})\nabla
\nabla\widetilde{g}(\mathbf{r}_{ij},\mathbf{r}^{\prime})\phi_{im}%
(\mathbf{r}^{\prime})\right]  \cdot\mathbf{c}_{im}\nonumber\\
&  +\frac{\omega^{2}\mu}{k^{2}}\sum_{m=1}^{M}{\int_{V_{\delta_{ij}}}%
\mathrm{d}\mathbf{r}^{\prime}\nabla^{2}g_{0}(\mathbf{r}_{ij},\mathbf{r}%
^{\prime})\left[  \Delta\epsilon(\mathbf{r}^{\prime})\phi_{im}(\mathbf{r}%
^{\prime})-\Delta\epsilon_{ij}\phi_{im}(\mathbf{r}_{ij})\right]  }%
\cdot\mathbf{c}_{im}.
\end{align}
}

When $n\neq i$, we have $V_{\delta_{ij}}\notin\Omega$, then the first integral
of (\ref{eqn:discre})
\begin{equation}
\int_{\Omega_{n}\backslash V_{\delta_{ij}}}\mathrm{d}\mathbf{r}^{\prime}%
\Delta\epsilon(\mathbf{r}^{\prime})\bar{\mathbf{G}}_{\mathbf{E}}%
(\mathbf{r}_{ij},\mathbf{r}^{\prime})\phi_{nm}(\mathbf{r}^{\prime})
\end{equation}
is regular in the whole domain $\Omega_{n}$ and hence it can be evaluated by
regular Gauss quadratures, i.e.,
\begin{equation}
\int_{\Omega_{n}}\mathrm{d}\mathbf{r}^{\prime}\Delta\epsilon(\mathbf{r}%
^{\prime})\bar{\mathbf{G}}_{\mathbf{E}}(\mathbf{r}_{ij},\mathbf{r}^{\prime
})\phi_{nm}(\mathbf{r}^{\prime})=\left(  \frac{a_{n}}{2}\right)  ^{3}%
\sum_{m=1}^{M}\Delta\epsilon_{nm}\bar{\mathbf{G}}_{\mathbf{E}}(\mathbf{r}%
_{ij},\mathbf{r}_{nm})\omega_{m}^{s},
\end{equation}
with $\omega_{m}^{s}$ being the standard Gauss weights in 3-D, which are
obtained from the tensor product of the Gauss weights in $[-1,1]$ or uniform
weights for the periodic direction in $\phi\in$ $[0,2\pi]$ in the cases of
spheres.
%\begin{eqnarray}\nonumber
%	\sum_{n=1}^N\sum_{m=1}^M\left[\int_{\Omega_n\backslash V_{\delta_{\tilde{i}\tilde{j}}}}\mathrm{d}{\bf r}'\Delta\epsilon\bar{\bf G}_{\bf E}({\bf r}_{\tilde{i}\tilde{j}}, {\bf r}')\phi_{ij}({\bf r}')\right]{\bf c}_{ij}&=&\sum_{j=1}^M\left[\int_{\Omega_{\tilde{i}}\backslash V_{\delta_{\tilde{i}\tilde{j}}}}\mathrm{d}{\bf r}'\Delta\epsilon\bar{\bf G}_{\bf E}({\bf r}_{\tilde{i}\tilde{j}}, {\bf r}')\phi_{\tilde{i}j}({\bf r}')\right]{\bf c}_{\tilde{i}j}\\\nonumber
%&+& \sum_{i=1, i\neq\tilde{i}}^N\sum_{j=1}^M\left[\int_{\Omega_i}\mathrm{d}{\bf r}'\Delta\epsilon\bar{\bf G}_{\bf E}({\bf r}_{\tilde{i}\tilde{j}}, {\bf r}')\phi_{ij}({\bf r}')\right]{\bf c}_{ij}
%\end{eqnarray}
%where $\Omega_i$ is the $i$-th element.

When $n=i$, although the singularity $\mathbf{r}_{ij}$ is excluded from the
domain $\Omega_{i}$, the calculation of the integral
\begin{equation}
\int_{\Omega_{i}\backslash V_{\delta_{ij}}}\mathrm{d}\mathbf{r}^{\prime}%
\Delta\epsilon(\mathbf{r}^{\prime})\bar{\mathbf{G}}_{\mathbf{E}}%
(\mathbf{r}_{ij},\mathbf{r}^{\prime})\phi_{im}(\mathbf{r}^{\prime})
\label{eqn:singular}%
\end{equation}
is still challenging. Therefore, we present an efficient quadrature formula to
evaluate this integral in the following subsection.

\subsection{Interpolated weights on tensor-product nodes for integrals on
$\Omega\backslash V_{\delta}$}

For the sake of generality, we consider the following integral:
\begin{equation}
I_{s}=\int_{\Omega\backslash V_{\delta}}\frac{f(\mathbf{r};\mathbf{r}^{\prime
})h(\mathbf{r};\mathbf{r}^{\prime})}{R^{k}}\mathrm{d}\mathbf{r}^{\prime}%
,\quad\mathbf{r}\in V_{\delta},\label{eqn:singular-general}%
\end{equation}
where $k=1,2,3$ corresponds to weak-, strong-, and hyper-singularity of the
integral, respectively. The function $f(\mathbf{r};\mathbf{r}^{\prime})$ is
assumed to be a general smooth and well-defined function, while $h(\mathbf{r}%
;\mathbf{r}^{\prime})$ is some fixed function resulting from the directional
derivative in the definition of the dyadic Green's function.

Because the function $f(\mathbf{r};\mathbf{r}^{\prime})$ is smooth over the
whole domain $\Omega$, it can be well approximated by the following simple
interpolation:
\begin{equation}
f(\mathbf{r};\mathbf{r}^{\prime})\approx\sum_{j=1}^{J}f(\mathbf{r}%
;\mathbf{r}_{j})\phi_{j}(\mathbf{r}^{\prime}),\text{ \ }\mathbf{r}_{j}%
\in\Omega, \label{eqn:interpolation}%
\end{equation}
where $\{\mathbf{r}_{j}\}_{j=1}^{J}$ are $J$-nodes in $\Omega$ formed by a
tensor product of Gauss nodes in $[-1,1]$ or uniformly spaced nodes in
$[0,2\pi]$ for periodic direction. Replacing $f(\mathbf{r};\mathbf{r}^{\prime
})$ in (\ref{eqn:singular-general}) with (\ref{eqn:interpolation}) yields
\begin{equation}
\int_{\Omega\backslash V_{\delta}}\frac{f(\mathbf{r};\mathbf{r}^{\prime
})h(\mathbf{r};\mathbf{r}^{\prime})}{R^{k}}d\mathbf{r}^{\prime}\approx
\sum_{j=1}^{J}f(\mathbf{r};\mathbf{r}_{j})\omega_{j},
\end{equation}
We call $\omega_{j}$ the interpolated weights, or just \emph{interpolated
weights}, defined through the integral
\begin{equation}
\omega_{j}=\int_{\Omega\backslash V_{\delta}}\frac{\phi_{j}(\mathbf{r}%
^{\prime})h(\mathbf{r};\mathbf{r}^{\prime})}{|\mathbf{r}-\mathbf{r}^{\prime
}|^{k}}d\mathbf{r}^{\prime}. \label{eqn:weights}%
\end{equation}
Note that the interpolated weights $\omega_{j}$ depend on the location of the
singularity $\mathbf{r}$ and relies on accurate calculations of
(\ref{eqn:weights}), which can be accomplished in two steps: for the case of a
cube domain $\Omega$, first $\Omega$ is subdivided upto 18 subcubes (the
number of subcubes depends on the location of the singularity) with one
containing the singularity $\mathbf{r}$ in its center. For the cube including
the singularity, a straight-forward, brute-force approach involving a large
number $N_{b}$ ($N_{b}\gg J$) of Gauss points is adopted in local spherical
polar coordinates to obtain satisfactory accuracy. While for the other cubes,
regular tensor product Gauss quadrature is applied. Details of the
computations and resulting weight $\left\{  \omega_{j}\right\}  $ tables and
also similar approaches for sphere and cylinder domains could be found in
\cite{Zinser:2015}.

The computation of weights $\omega_{j}$ only needs to be performed once and
tabulated for the reference domain; they can be used for general cubic or
spherical domains. Due to the smoothness of function $f(\mathbf{r}%
,\mathbf{r}^{\prime})$, the number $J$ is relatively small, especially if the
element size is small as in meta-material designs, so the computation of
(\ref{eqn:singular}) is efficient once $\omega_{j}$ are obtained.

\begin{rem}
As the CPV is used in the VIE, the VIE depends on the specific shape of the
exclusion volume. Thus, the interpolated weights can only be used for the
shape for which it was calculated, namely, a cube or a sphere shape here. So,
for a general element obtained by an affine mapping as in a finite element
triangulation, we will need to isolate the singularity by a cube or a sphere
with the singularity at its center, then the pre-calculated interpolated
weights defined in (\ref{eqn:weights}) can be used where integral over the
rest of the region within the element of a more general shape can be computed
with regular Gauss quadratures.
\end{rem}

\subsection{Computation of VIE matrix entries}

In this section, we will show how to compute the matrix entries accurately for
the VIE in the following steps.

\begin{itemize}
\item Step I: calculate interpolated weights $\left\{  \omega_{j}\right\}  $
on the reference (cubic or spherical) domain.

\bigskip

For Eq. (\ref{eqn:weights}), we take $V_{\delta}=B(\mathbf{r}_{j},{\delta})$,
where $\delta>0$ is a prescribed small quantity and $\mathbf{r}_{j}%
,j=1,2,...M$ are the $M$ tensor-product nodes in the reference domain, and we
let $J=M$ in Eq. (\ref{eqn:interpolation}).

Next, we have the dyadic Green's function for the free-space in the following
form,
\begin{align}
\overline{\mathbf{G}}_{\mathbf{E}}  &  =g\mathbf{I}+\frac{\nabla^{2}g}{k^{2}%
}=\frac{e^{-ikR}}{4\pi R}(\mathbf{I}-\mathbf{u}\otimes\mathbf{u}%
)\nonumber\label{eqn:G-full}\\
&  -\frac{ie^{-ikR}}{4\pi R^{2}k}(\mathbf{I}-3\mathbf{u}\otimes\mathbf{u}%
)-\frac{e^{-ikR}}{4\pi R^{3}k^{2}}(\mathbf{I}-3\mathbf{u}\otimes\mathbf{u}).
\end{align}
%where $R=|{\bf r}-{\bf r}'|$ and ${\bf u}$ is the unit vector of ${\bf r}-{\bf r}'$, i.e., ${\bf u}=\displaystyle{\frac{{\bf r}-{\bf r}'}{R}}$,
Since the value of function $\mathbf{u}\otimes\mathbf{u}$ is multi-valued at
$R=0$, in (\ref{eqn:singular-general}) we will take
\begin{equation}
h(\mathbf{r};\mathbf{r}^{\prime})=\mathbf{u}\otimes\mathbf{u,}%
\end{equation}
and hence (\ref{eqn:weights}) will produce a set of $9$ interpolated weights.
However, due to the symmetry of the matrix $\mathbf{u}\otimes\mathbf{u}$, only
6 components need to be considered. For the identity matrix $\mathbf{I}$ term
in (\ref{eqn:G-full}), we need a set of scalar interpolation weights by
assuming $h(\mathbf{r},\mathbf{r}^{\prime})=1$ in (\ref{eqn:weights}).

Additionally, for the scalar and matrix weights, we will consider $k=1,2$, and
$3$ for weak, strong, and hyper singular integrals, respectively.

In summary, for each collocation point (also the singularity location)
$\mathbf{r}_{j},j=1,2,...,M$ in an element, scalar weights $\omega_{j,m}$,
$\bar{\omega}_{j,m}$ and $\tilde{\omega}_{j,m},m=1,2,...,M$ are calculated for
weak-, strong-, hyper-singular integrals, respectively. And the corresponding
matrix weights are denoted as $\Lambda_{j,m}$, $\bar{\Lambda}_{j,m}$, and
$\tilde{\Lambda}_{j,m}$. Here the first index $j$ indicates the location of
the singularity of the integrand while the second one $m$ is the quadrature
weight index. These weights only need to be calculated once and then stored
for future use.

\bigskip

\item Step II: We consider the computational domain consisting of a collection
of fundamental elements with size $a_{i},i=1,2,...,N$ and assign $M$
collocation points in each elements. For the $j$-th collocation point
$\mathbf{r}_{ij}$ in the $i$-th element, we construct the equation:
\begin{equation}
-\omega^{2}\mu\sum_{n=1}^{N}\sum_{m=1}^{M}\mathbf{A}_{nm}\cdot\mathbf{c}%
_{nm}-\sum_{m=1}^{M}\mathbf{B}_{im}\cdot\mathbf{c}_{im}+\left(  1+\frac{1}%
{3}\Delta\epsilon_{ij}\right)  \mathbf{I}_{3\times3}\cdot\mathbf{c}%
_{ij}=\mathbf{E}_{ij}^{\mathrm{inc}}, \label{eqn:equation}%
\end{equation}
%		The coefficients
%		\begin{equation}
%			{\bf C}=\left(1+\frac{1}{3}\Delta\epsilon_{ij}\right){\bf I}
%		\end{equation}
The matrix $\mathbf{B}$ originates from the correction terms of the Cauchy
principal value
\begin{align}
\mathbf{B}_{im}  &  =\omega^{2}\mu\int_{B(\mathbf{r}_{ij},a_{i}\delta
)}\mathrm{d}\mathbf{r}^{\prime}\Delta\epsilon(\mathbf{r}^{\prime}%
)g(\mathbf{r}_{ij},\mathbf{r}^{\prime})\phi_{im}(\mathbf{r}^{\prime
})\nonumber\\
&  +\frac{\omega^{2}\mu}{k^{2}}\int_{B(\mathbf{r}_{ij},a_{i}\delta)}%
\mathrm{d}\mathbf{r}^{\prime}\Delta\epsilon(\mathbf{r}^{\prime})\nabla
^{2}\tilde{g}(\mathbf{r}_{ij},\mathbf{r}^{\prime})\phi_{im}(\mathbf{r}%
^{\prime})\nonumber\\
&  +\frac{\omega^{2}\mu}{k^{2}}\int_{B(\mathbf{r}_{ij},a_{i}\delta)}%
\mathrm{d}\mathbf{r}^{\prime}\nabla^{2}g_{0}(\mathbf{r}_{ij},\mathbf{r}%
^{\prime})\left[  \Delta\epsilon(\mathbf{r}^{\prime})\phi_{im}(\mathbf{r}%
^{\prime})-\Delta\epsilon_{ij}\phi_{im}(\mathbf{r}_{ij}))\right]  ,
\end{align}
and it can be calculated by standard Gauss quadrature through spherical
coordinates since the Jacobian will eliminate completely the singularity of
the integrands.

When $n=i$, we calculate the integral (\ref{eqn:singular}) as
\begin{align}
\mathbf{A}_{im}  &  =\frac{1}{4\pi}J_{i}\sum_{j=1}^{M}\Delta\epsilon
_{im}\left[  \left(  e^{-ikR_{m}}\omega_{j,m}^{i}-ie^{-ikR_{m}}\bar{\omega
}_{j,m}^{i}-e^{-ikR_{m}}\tilde{\omega}_{j,m}^{i}\right)  \mathbf{I}_{3\times
3}\right. \nonumber\\
&  \left.  e^{-ikR_{m}}\Lambda_{j,m}^{i}-ie^{-ikR_{m}}\bar{\Lambda}_{j,m}%
^{i}-e^{-ikR_{m}}\tilde{\Lambda}_{j,m}^{i}\right]  ,
\end{align}
where $R_{m}=|\mathbf{r}_{ij}-\mathbf{r}_{im}|$ and $J_{i}$ is the Jacobian
from the reference domain to the physical domain $\Omega_{i}$.

Note that the interpolated quadrature weights are rescaled from the reference
domain. In the cube example, the reference domain is $[-1,1]^{3}$ and if the
physical domain has length $a_{i}$, then $J_{i}=\displaystyle{\frac{a_{i}}{2}%
}$ and recall the definition of the interpolated weights (\ref{eqn:weights}),
we have
\begin{equation}%
\begin{array}
[c]{lll}%
\omega_{j,m}^{i}=\left(  \frac{2}{a_{i}}\right)  \omega_{j,m}, & \bar{\omega
}_{j,m}^{i}=\left(  \frac{2}{a_{i}}\right)  ^{2}\bar{\omega}_{j,m}, &
\tilde{\omega}_{j,m}^{i}=\left(  \frac{2}{a_{i}}\right)  ^{3}\tilde{\omega
}_{j,m}\\
\label{eqn:interpolated-weights}\Lambda_{m}^{i}=\left(  \frac{2}{a_{i}%
}\right)  \Lambda_{j,m}, & \bar{\Lambda}_{j,m}^{i}=\left(  \frac{2}{a_{i}%
}\right)  ^{2}\bar{\Lambda}_{j,m}, & \tilde{\Lambda}_{j,m}^{i}=\left(
\frac{2}{a_{i}}\right)  ^{3}\tilde{\Lambda}_{j,m}%
\end{array}
,
\end{equation}
when $n\neq i$, we have
\begin{equation}
\mathbf{A}_{nm}=J_{n}\sum_{j=1}^{M}\Delta\epsilon_{nm}\bar{\mathbf{G}%
}_{\mathbf{E}}(\mathbf{r}_{ij},\mathbf{r}_{nm})\omega_{j}^{s}.
\end{equation}

%Map the reference domain to the physical element with length $a$ excluding a ball with radius $a\delta$.

\item Step III: (\ref{eqn:equation}) for all the $N\times M$ tensor-product
nodes can be assembled as the following linear algebraic equation system
\begin{equation}
\mathbf{V}\cdot\vec{\mathbf{c}}=\left[
\begin{array}
[c]{ccc}%
\mathbf{V}_{xx} & \mathbf{V}_{xy} & \mathbf{V}_{xz}\\
\mathbf{V}_{yx} & \mathbf{V}_{yy} & \mathbf{V}_{yz}\\
\mathbf{V}_{zx} & \mathbf{V}_{zy} & \mathbf{V}_{zz}%
\end{array}
\right]  \cdot\left[
\begin{array}
[c]{c}%
\mathbf{c}_{x}\\
\mathbf{c}_{y}\\
\mathbf{c}_{z}%
\end{array}
\right]  =\left[
\begin{array}
[c]{c}%
\mathbf{E}_{x}^{\mathrm{inc}}\\
\mathbf{E}_{y}^{\mathrm{inc}}\\
\mathbf{E}_{z}^{\mathrm{inc}}%
\end{array}
\right]  . \label{eqn:sys}%
\end{equation}

\end{itemize}

Based on the properties of the Green's function, the $3NM\times3NM$ matrix
$\mathbf{V}$ is partitioned into nine blocks, each block is a $NM\times NM$
sub-matrix. The solution of the VIE contains three $NM\times1$ vectors, which
represent the field in $x$, $y$, and $z$ directions. Then, the system is
solved by a matrix solver. In this paper, we used the generalized minimal
residual (GMRES) method \cite{Saad}.

\section{Numerical Results}

\label{sec:results}

In this section we test the accuracy of the interpolated weights on
tensor-product nodes, the $\delta$-independence of the solution of the VIE,
and convergence of the $p$-refinement of the Nystr\"{o}m collocation method.

\subsection{Accuracy of the interpolated weights on tensor-product nodes for
computing matrix entries}

In this subsection, the accuracy of the interpolated weights on tensor-product
on a cube is presented. The study on a sphere can be treated in a similar way.

In (\ref{eqn:equation}), the calculation of matrix $\mathbf{B}$ from the
correction terms is straightforward; so we will focus on validating the
interpolated weights in computing matrix $\mathbf{A}$. For convenience, we
consider the integral of a real-valued, tensor function
\begin{equation}
\frac{\cos{R}}{R}(\mathbf{I}-\mathbf{u}\otimes\mathbf{u})+\frac{\cos{R}}%
{R^{2}}(\mathbf{I}-3\mathbf{u}\otimes\mathbf{u})+\frac{\cos{R}}{R^{3}%
}(\mathbf{I}-3\mathbf{u}\otimes\mathbf{u}), \label{eqn:sampleintegral}%
\end{equation}
which is similar to the Green's function in (\ref{eqn:G-full}) in the domain
$\Omega\backslash V_{\delta}$. Without loss of generality, we take
$\Omega=[-1,1]^{3}$ and $\mathbf{r}_{j}$ as the 27 points constructed from the
tensor product of the Gauss points of order 3 in $[-1,1]$. Thus, we have
$j=m=1,2,...,27$ as in (\ref{eqn:interpolated-weights}) and the integral can
be written as
\begin{equation}
\mathbf{G}_{j}\approx\sum_{m=1}^{27}\cos(|\mathbf{r}_{m}-\mathbf{r}%
_{j}|)\left[  \left(  \omega_{j,m}+\bar{\omega}_{j,m}+\tilde{\omega}%
_{j,m}\right)  \bar{\mathbf{I}}-\Lambda_{j,m}-3\bar{\Lambda}_{j,m}%
-3\tilde{\Lambda}_{j,m}\right]  , \label{eqn:appxofG}%
\end{equation}
which is a $3\times3$ matrix depending on $\mathbf{r}_{j}$.
%		we test the accuracy of the interpolated Gauss weights by the sample integral
%		\begin{equation}\label{eqn:sampleintegral}
%			{\bf G}=\int_{\Omega\backslash V_{\delta}}\left[\frac{\cos{R}}{R}({\bf I}-{\bf u}\otimes{\bf u})+\frac{\cos{R}}{R^2}({\bf I}-3{\bf u}\otimes{\bf u})+\frac{\cos{R}}{R^3}({\bf I}-3{\bf u}\otimes{\bf u})\right]\mathrm{d}{\bf r}',
%		\end{equation}
%		
%		  and Eq. (\ref{eqn:sampleintegral}) is approximated as

For each $\mathbf{G}_{j}$, we use the direct brute force method introduced in
\cite{Zinser:2015} to obtain the reference value of the matrix entries with a
small $\delta=10^{-3}$. Then we calculate the integral using the interpolated
weights as in (\ref{eqn:appxofG}) with different values of $\delta$. According
to the previous analysis, the error should decay on the order $O(\delta^{2})$.

We classify the 27 sets of weights into four categories, based on the position
of the singularity $\mathbf{r}_{j}$, as being near the corner, edge, face, and
center of the cube. The matrix $\mathbf{G}_{j}$ is symmetric, so we only check
the three diagonal entries ($g_{11}$, $g_{22}$, and $g_{33}$) and the three
upper triangular entries (${g_{12}}$, $g_{13}$, and $g_{23}$).

Table \ref{table:center} presents the numerical results when the singularity
$\mathbf{r}_{j}$ is located in the center of the cube, in which case
$g_{11}=g_{22}=g_{33}$ and the off-diagonal entries are all zeros. The
reference solution is $g_{11}=4.027477$ while the $g_{11}$ computed with
$\delta=$ $0.1,0.05,0.025$, and $0.0125$ are $3.985701,4.017024,4.024872,$and
$4.026835$, respectively.

\begin{table}[ptb]
\caption{Convergence of the integral when the singularity is at center.
$g_{11}=g_{22}=g_{33}$ and $g_{12}=g_{13}=g_{23}=0$.}%
\label{table:center}%
\centering
\begin{tabular}
[c]{l|l|l|l|l}\hline\hline
& $\delta= 0.1$ & $\delta= 0.05$ & $\delta= 0.025$ & $\delta= 0.0125$\\\hline
$g_{11}$ & 3.985701 & 4.017024 & 4.024872 & 4.026835\\
error & -4.1784E-2 & -1.0461E-2 & -2.613E-3 & -6.5E-4\\
order & - & 2 & 2 & 2\\\hline\hline
\end{tabular}
\end{table}

\begin{table}[ptb]
\caption{Convergence of the integral when singularity is at a corner.
Reference solution $g_{11}=g_{22}=g_{33}=0.982526$ and $g_{12}=g_{13}%
=g_{23}=-0.998097$.}%
\label{table:corner}%
\centering
\begin{tabular}
[c]{l|l|l|l|l}\hline\hline
& $\delta= 0.1$ & $\delta= 0.05$ & $\delta= 0.025$ & $\delta= 0.0125$\\\hline
$g_{11}$ & 0.940714 & 0.972063 & 0.979913 & 0.981876\\
error & -2.80425E-1 & -7.3424E-2 & -1.8102E-2 & -3.983E-3\\
order & - & 1.93 & 2 & 2.1\\\hline\hline
$g_{12}$ & -0.998084 & -0.998094 & -0.998097 & -0.998097\\
error & 1.3E-5 & 3.0E-6 & 0 & 0\\
order & - & 2.1 & - & -\\\hline\hline
\end{tabular}
\end{table}

\begin{table}[ptb]
\caption{Convergence of the integral when singularity is at an edge. Reference
solutions $g_{11}=-1.515302$, $g_{22}=g_{33}=3.39234$, $g_{23}=-1.579086$ and
$g_{12}=g_{13}=0$.}%
\label{table:edge}%
\centering
\begin{tabular}
[c]{l|l|l|l|l}\hline\hline
& $\delta= 0.1$ & $\delta= 0.05$ & $\delta= 0.025$ & $\delta= 0.0125$\\\hline
$g_{11}$ & -1.559532 & -1.526351 & -1.518059 & -1.515987\\
error & -4.423E-2 & 1.1049E-2 & 2.757E-3 & -6.85E-4\\
order & - & 2 & 2 & 2\\\hline\hline
$g_{22}$ & 3.35175 & 3.38217 & 3.389798 & 3.391707\\
error & -4.059E-2 & -1.1017E-2 & 2.542E-3 & -6.33E-4\\
order & - & 1.9 & 2.1 & 2\\\hline\hline
$g_{23}$ & -1.579072 & -1.579082 & -1.579085 & -1.579085\\
error & 1.4E-5 & 4.0E-6 & 1.0E-6 & 1.0E-6\\
order & - & 1.8 & 2 & 0\\\hline\hline
\end{tabular}
\end{table}

\begin{table}[ptb]
\caption{Convergence of the integral when singularity is at a face. Reference
solutions $g_{11}=0.877428$, $g_{22}=0$, $g_{33}=6.494784$, and $g_{12}%
=g_{13}=g_{23}=0$.}%
\label{table:face}%
\centering
\begin{tabular}
[c]{l|l|l|l|l}\hline\hline
& $\delta= 0.1$ & $\delta= 0.05$ & $\delta= 0.025$ & $\delta= 0.0125$\\\hline
$g_{11}$ & 0.83442 & 0.866672 & 0.874742 & 0.87676\\
error & -4.3008E-2 & -1.10756E-2 & -2.686E-3 & -6.68E-4\\
order & - & 1.9 & 2 & 2\\\hline\hline
$g_{33}$ & 6.455419 & 6.484909 & 6.492315 & 6.49417\\
error & -3.9365E-2 & -9.875E-3 & -2.469E-3 & -6.14E-4\\
order & - & 2 & 2 & 2\\\hline\hline
\end{tabular}
\end{table}

In a similar fashion, Tables \ref{table:corner}$-$\ref{table:face} show the
accuracy when the singularity is located near the corner, edge, and face of
the cube, respectively. When compared to the reference values, the expected
$O(\delta^{2})$ behavior is confirmed.

\bigskip

\subsection{Exclusion volume $\delta$-independence of the VIE solution}

\begin{figure}[ptb]
\begin{center}%
\begin{tabular}
[c]{cc}%
\includegraphics[width=0.5\textwidth]{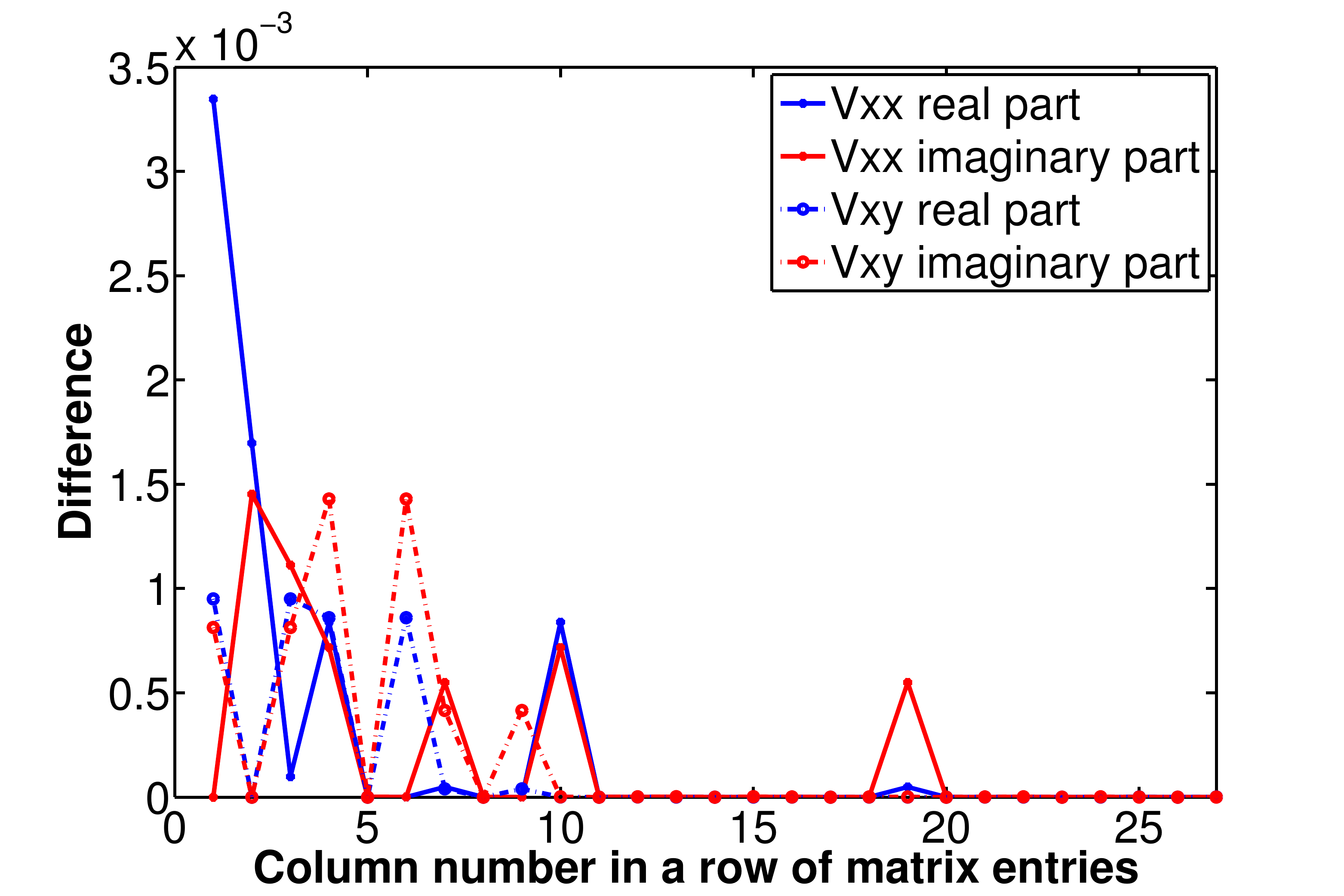} &
\includegraphics[width=0.5\textwidth]{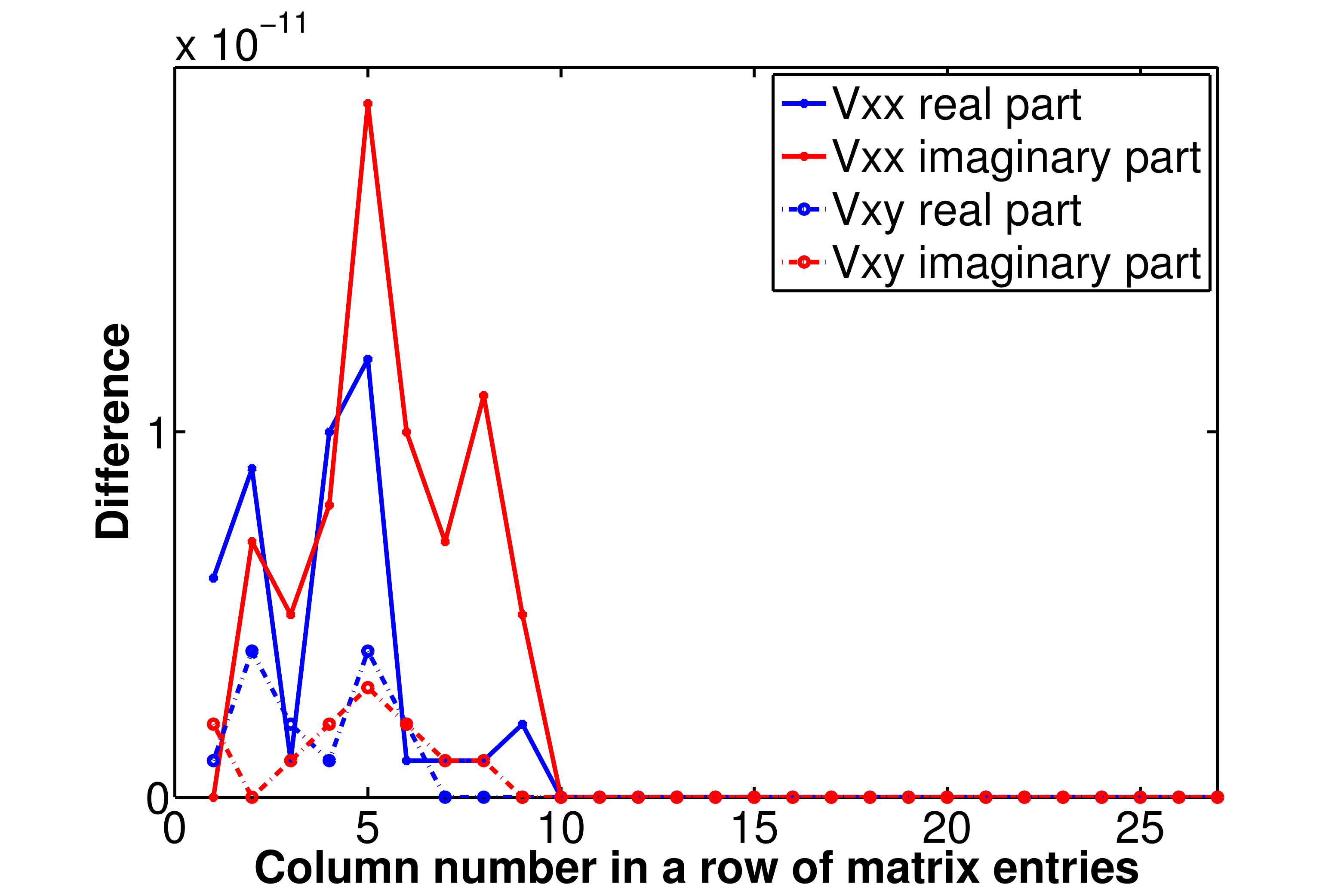}\\
(a) & (b)
\end{tabular}
\end{center}
\caption{Differences of matrix entries with $\delta$=0.1 and $\delta=0.001$.
(a): without correction terms; (b) with correction terms.}%
\label{fig:entryerror}%
\end{figure}

Equation (\ref{eqn:per-217}) provides a formulation with which the solution of
the VIE will be independent of the choice of the exclusion volume size
$\delta$. In the following tests, we take $\mu=1$, $\Delta\epsilon=4$,
$\omega=1$ and solve the VIE. The computational domain is taken as
$[-\pi/2,\pi/2]^{3}$, while the incident wave is
\begin{equation}
\mathbf{E}_{x}^{\mathrm{inc}}=e^{ik(-y+0.5z)},\quad\mathbf{E}_{y}%
^{\mathrm{inc}}=\mathbf{E}_{z}^{\mathrm{inc}}=0.\nonumber
\end{equation}
We first check the $\delta$-independence of the matrix entries in
(\ref{eqn:sys}). Fig. \ref{fig:entryerror} displays the differences of one row
of entries in the matrix $\mathbf{V}$ between the choices of $\delta=0.1$ and
$\delta=0.001$, in which the solid lines are for the entries from a diagonal
block ($V_{xx}$) and dashed lines are for the entries from an off-diagonal
block ($V_{xy}$). The blue curves are for real while red curves are for
imaginary parts. From Fig. \ref{fig:entryerror}(a) we can see that the
differences between entries in the corresponding positions can be as large as
$3.0\times10^{-3}$ when the correction terms are not included. In contrast,
the corresponding differences are reduced to below $2\times10^{-11}$ when the
corrections are included. Hence, the matrix entries are $\delta$-independent
when the correction terms are included. 

Next we check the $\delta$-dependence of the overall solution of the VIE. The
solution of VIE with a very small $\delta=0.001$ is chosen as the reference
solution. Then, the numerical solutions are computed with $\delta$
$=0.1,0.05,0.025,$ and $0.0125$ and compared with the reference solution. The
differences are measured in the $L^{\infty}$ norm for the three components
$E_{x}$, $E_{y}$, $E_{z}$ and they are listed without and with the correction
terms in Tables \ref{table:sol1} and \ref{table:sol2}, respectively.
\begin{table}[th]
\caption{Comparison of solutions of the VIE without the correction terms}%
\label{table:sol1}%
\centering
\begin{tabular}
[c]{l|l|l|l|l}\hline\hline
& $\delta= 0.1$ & $\delta= 0.05$ & $\delta= 0.025$ & $\delta= 0.0125$\\\hline
$\|E_{x}-E_{x}^{\mathrm{ref}}\|_{L^{\infty}}$ & 3.360E-3 & 8.264E-4 &
2.044E-4 & 5.039E-5\\
$\|E_{y}-E_{y}^{\mathrm{ref}}\|_{L^{\infty}}$ & 1.476E-3 & 3.696E-4 &
9.258E-5 & 2.358E-5\\
$\|E_{z}-E_{z}^{\mathrm{ref}}\|_{L^{\infty}}$ & 2.533E-3 & 6.387E-4 &
1.597E-4 & 3.969E-5\\\hline\hline
\end{tabular}
\end{table}\begin{table}[thpth]
\caption{Comparison of solutions of the VIE with the correction terms}%
\label{table:sol2}%
\centering
\begin{tabular}
[c]{l|l|l|l|l}\hline\hline
& $\delta= 0.1$ & $\delta= 0.05$ & $\delta= 0.025$ & $\delta= 0.0125$\\\hline
$\|E_{x}-E_{x}^{\mathrm{ref}}\|_{L^{\infty}}$ & 8.0E-12 & 1.0E-12 & 0 & 0\\
$\|E_{y}-E_{y}^{\mathrm{ref}}\|_{L^{\infty}}$ & 2.0E-12 & 1.0E-12 & 0 & 0\\
$\|E_{z}-E_{z}^{\mathrm{ref}}\|_{L^{\infty}}$ & 1.0E-12 & 0 & 0 &
0\\\hline\hline
\end{tabular}
\end{table}

From Table \ref{table:sol1} it can be seen that without the correction terms,
the solution of VIE has an obvious dependence on the choice of $\delta$ and
the differences decreases in the order of $O(\delta^{2})$, while the solution
is indeed $\delta$-independent when the correction terms are included, as
shown in Table \ref{table:sol2}.

\subsection{$p$-convergence of the VIE solution}

For numerical convergence study, we can either discretize the computational
domain into smaller elements (increasing the number $N,$ $h$-refinement) or
use higher order polynomial basis (increasing the number $p,$ $p$-refinement)
in computing the integral in the VIE. In the current work we focus on the
latter since the emphasis is on the accurate calculation of the CPV of the
integral. We will check the solution of VIE on a sphere and a cube. To verify
the rate of convergence, we consider the linear relation between the
$\log_{10}$ of the energy error and $p$, the effective number of collocation
points in one direction, i.e.
\begin{equation}
\log_{10}(\text{Error})=\alpha p+\beta,\nonumber
\end{equation}
where
\begin{equation}
\text{Error}=\Vert\mathbf{E}^{p}-\mathbf{E}^{\mathrm{ref}}\Vert_{L^{2}%
(\Omega)},\nonumber
\end{equation}
and $\mathbf{E}^{p}$ is the solution computed with $p$ collocation points in
one direction, $\mathbf{E}^{\mathrm{ref}}$ is the reference solution, and
$\alpha$ and $\beta$ are fitting parameters. In previous sections we know that
the solution does not depend on the choice of parameter $\delta$, so we take
$\delta=0.001$ for the rest of calculations.

\begin{itemize}
\item Case 1: Solution of the VIE for the MIE scattering of a sphere
\end{itemize}

%\noindent\ :

The analytic solution of the VIE in a sphere is given by the Mie series
\cite{Mie}. For an incident wave
\begin{equation}
\mathbf{E}^{\mathrm{inc}}=\mathbf{i}_{x}e^{ikz},
\end{equation}
where $\mathbf{i}_{x}$ is the unit vector along $x$-direction in the Cartesian
coordinates, the exact solution of electromagnetic fields inside a sphere can
be expressed as the following series, \emph{in unit vectors of spherical
coordinates},
\begin{equation}
\mathbf{E}(r,\theta,\phi)=\sum_{n=1}^{\infty}\frac{i^{n}(2n+1)}{n(n+1)}\left(
c_{n}\mathbf{M}_{o1n}^{(1)}-id_{n}\mathbf{N}_{e1n}^{(1)}\right)  ,
\label{eqn:Mie}%
\end{equation}
where the coefficients are
\begin{equation}%
\begin{array}
[c]{l}%
c_{n}=\displaystyle{\frac{j_{n}(ka)[kah_{n}^{(1)}(ka)]^{\prime}-h_{n}%
^{(1)}(ka)[kaj_{n}(ka)]^{\prime}}{j_{n}(mka)[kah_{n}^{(1)}(ka)]^{\prime}%
-h_{n}^{(1)}(ka)[mkaj_{n}(mka)]^{\prime}}};\\
d_{n}=\displaystyle{\frac{mj_{n}(ka)[kah_{n}^{(1)}(ka)]^{\prime}-mh_{n}%
^{(1)}(ka)[kaj_{n}(ka)]^{\prime}}{m^{2}j_{n}(mka)[kah_{n}^{(1)}(ka)]^{\prime
}-h_{n}^{(1)}(x)[mkaj_{n}(mka)]^{\prime}}}.
\end{array}
\end{equation}
and the vector special harmonics are defined by
\begin{equation}
\mathbf{{M}_{o1n}^{(1)}=\left(
\begin{array}
[c]{c}%
0\\
\cos{\phi}\cdot\pi_{n}(\cos{\theta})j_{n}(mkr)\\
-\sin{\phi}\cdot\tau_{n}(\cos{\theta})j_{n}(mkr)
\end{array}
\right)  ,}%
\end{equation}

\begin{equation}
\mathbf{{N}_{e1n}^{(1)}=\left(
\begin{array}
[c]{c}%
n(n+1)\cos{\phi}\cdot\sin{\theta}\cdot\pi_{n}(\cos{\theta})\displaystyle{\frac
{j_{n}(mkr)}{mkr}}\\
\cos{\phi}\cdot\tau_{n}(\cos{\theta})\displaystyle{\frac{[mkrj_{n}%
(mkr)]^{\prime}}{mkr}}\\
-\sin{\phi}\cdot\pi_{n}(\cos{\theta})\displaystyle{\frac{[mkrj_{n}%
(mkr)]^{\prime}}{mkr}}%
\end{array}
\right)  .}%
\end{equation}
In the above formulas, $m$ is the refractive index of the sphere relative to
the ambient medium, $a$ the radius of the sphere and $k$ is the wave number of
the ambient medium. The functions $j_{n}(z)$ and $h_{n}^{(1)}(z)$ are
spherical Bessel functions of first and third kind, respectively, and their
derivatives have the relations \cite{Abramowitz}
\begin{equation}
\lbrack zj_{n}(z)]^{\prime}=zj_{n-1}(z)-nj_{n}(z);[zh_{n}^{(1)}(z)]^{\prime
}=zh_{n-1}^{(1)}(z)-nh_{n}^{(1)}(z).
\end{equation}
$\pi_{n}(\cos{\theta})$ and $\tau_{n}(\cos{\theta})$ have the relations
\begin{equation}
\pi_{n}=\frac{2n-1}{n-1}\cos{\theta}\cdot\pi_{n-1}-\frac{n}{n-1}\pi_{n-2}%
;\tau_{n}=n\cos{\theta}\cdot\pi_{n}-(n+1)\pi_{n-1},
\end{equation}
with
\begin{equation}
\pi_{0}=0;\pi_{1}=1;\pi_{2}=3\cos{\theta};\tau_{0}=0;\tau_{1}=\cos{\theta
};\tau_{2}=3\cos{2\theta}.
\end{equation}

To find the interpolated quadrature weights in (\ref{eqn:weights}), Lagrange
interpolation is used along the $r$-direction and $m_{r}$ Gauss nodes are used
while Fourier interpolation is applied for the $\phi\in\lbrack0,2\pi)$ and
$\theta\in\lbrack0,\pi]$. In total $m_{\theta}+1$ and $2m_{\phi}$ grid points
(namely, the collocation points in the Nystr\"{o}m method) are equally
distributed for $\theta$ and $\phi$, respectively. Therefore, the equivalent
number of collocation points in each direction is calculated as $p=\sqrt[3]%
{m_{r}(2m_{\phi}(m_{\theta}-1)+2)}$. The specific quadrature weights for these
nodes can be found in \cite{Zinser:2015}. \begin{figure}[ptb]
\begin{center}
\includegraphics[width=0.6\textwidth]{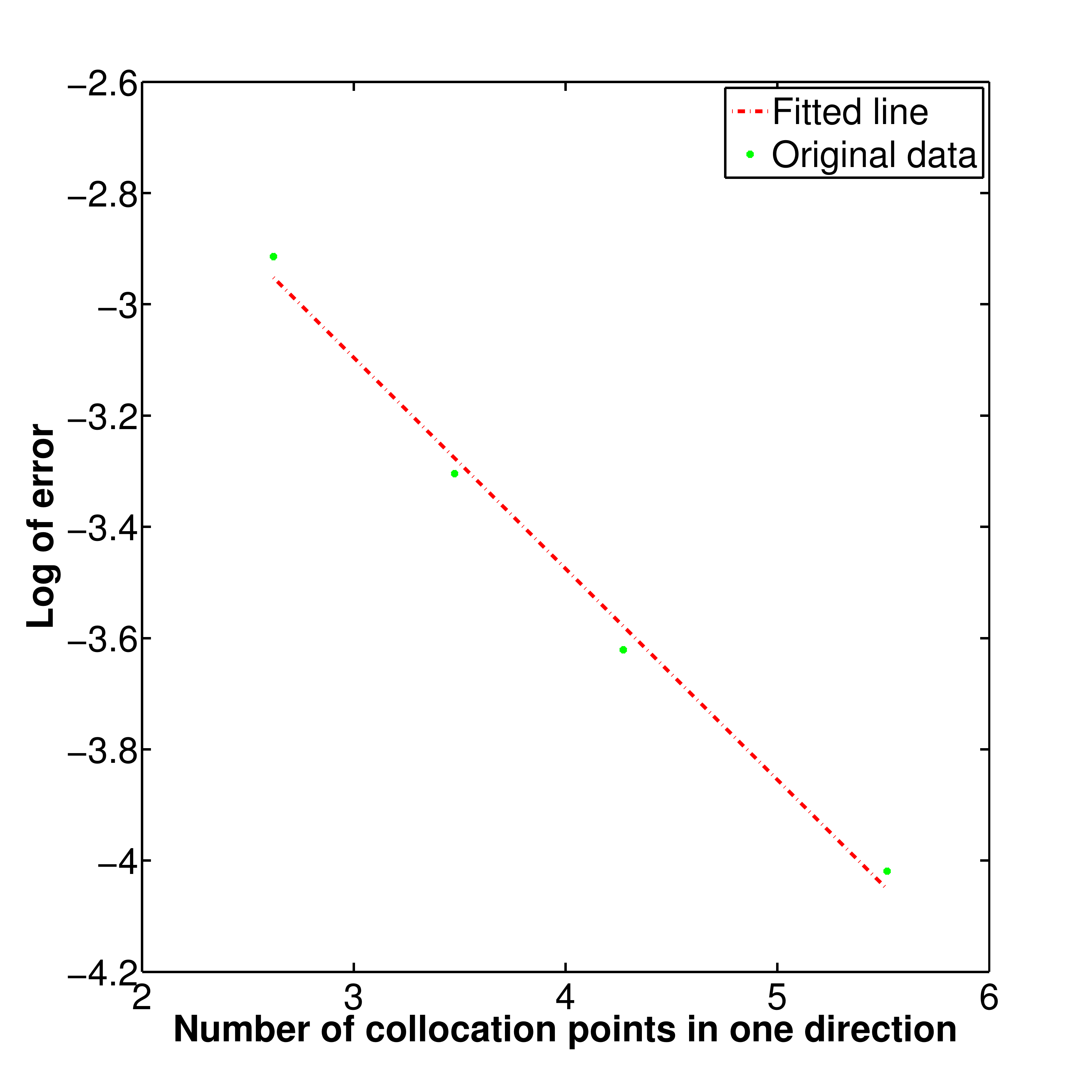}
\end{center}
\caption{$\log_{10}$ errors of solutions against number of collocation points
along one direction for $\mathbf{E}$ in a sphere.}%
\label{fig:p-conv-sphere}%
\end{figure}

\begin{figure}[ptb]
\begin{center}%
\begin{tabular}
[c]{cccc}%
\includegraphics[width=0.24\textwidth]{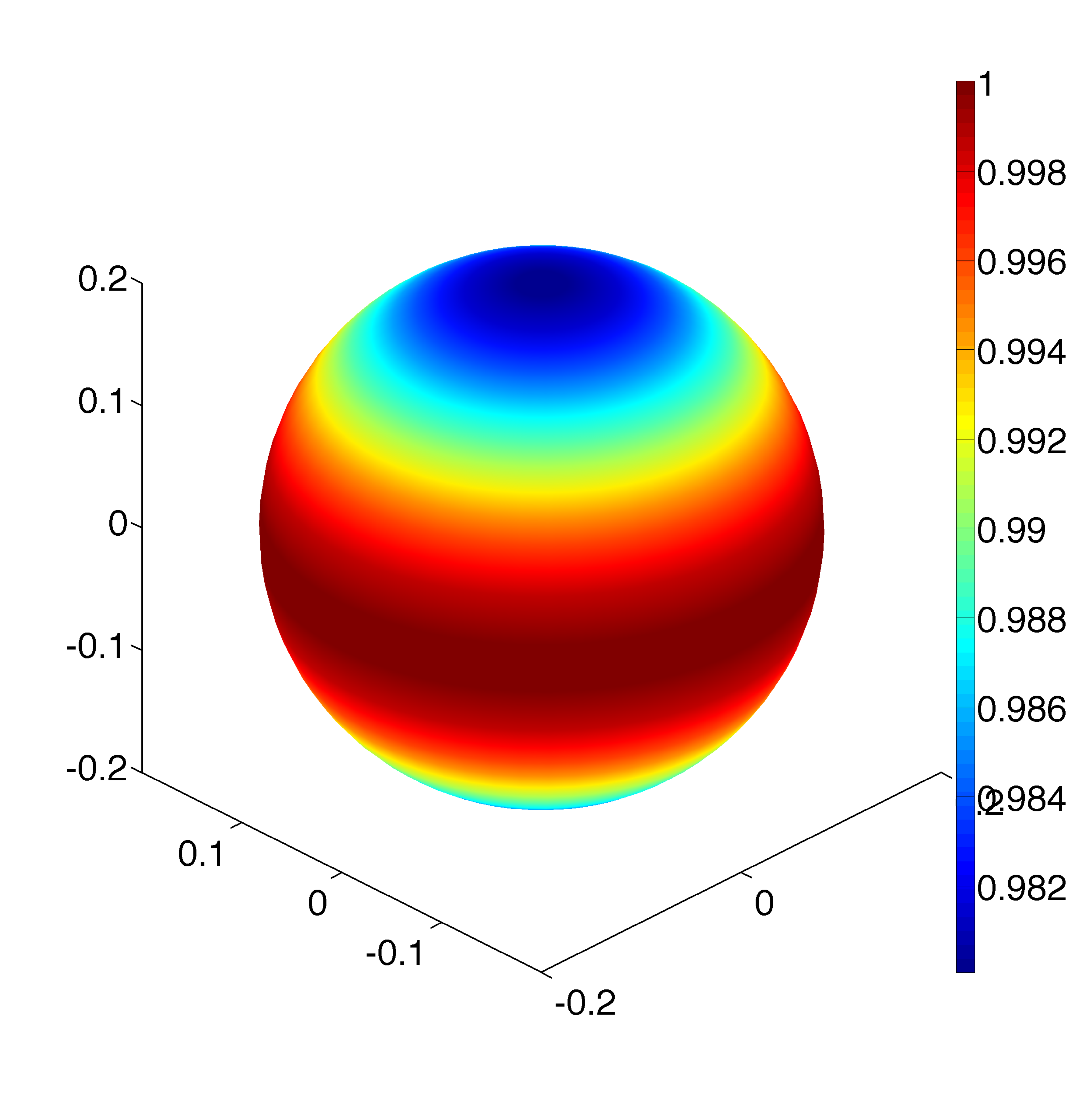} &
\includegraphics[width=0.24\textwidth]{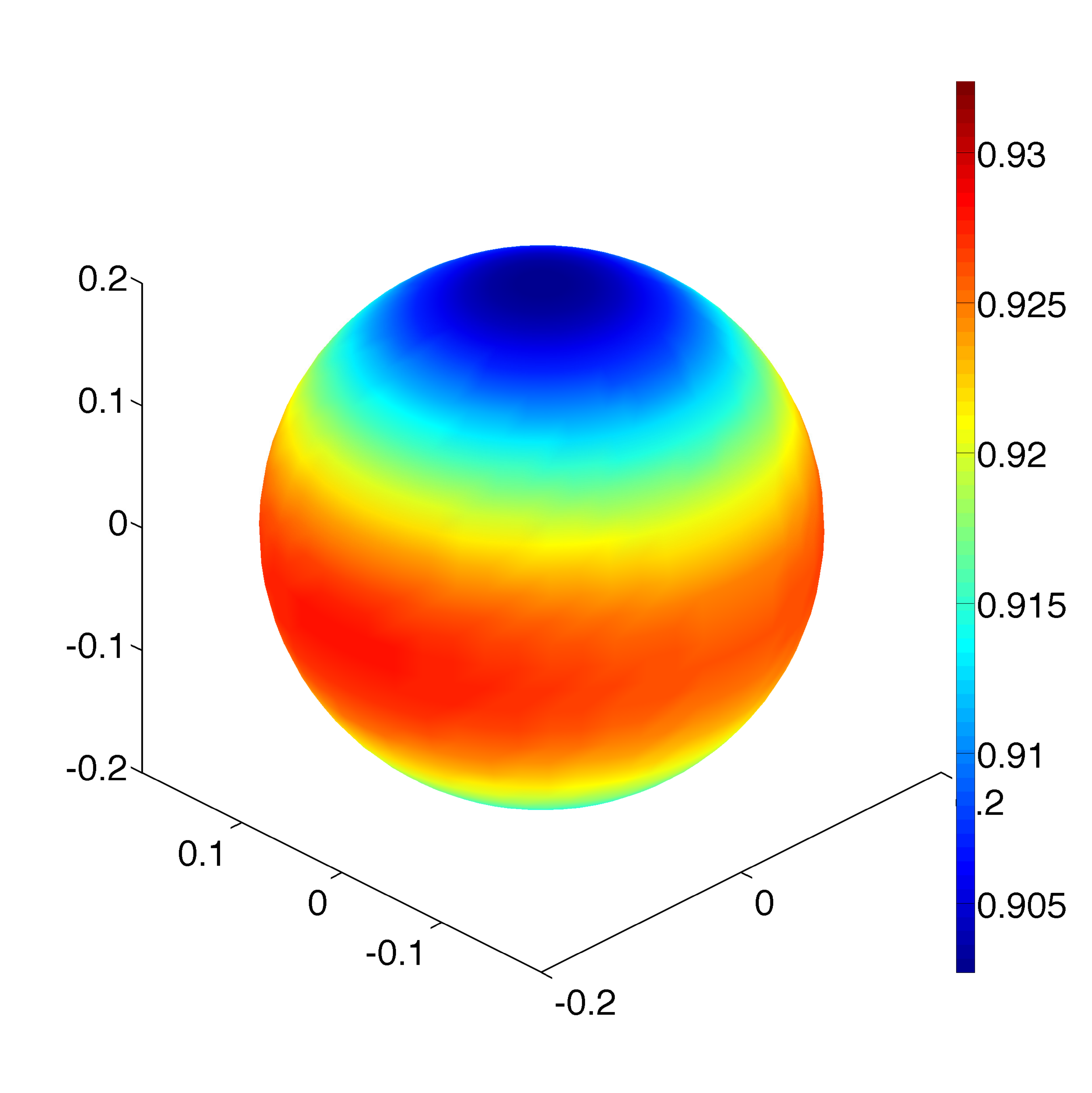} &
\includegraphics[width=0.24\textwidth]{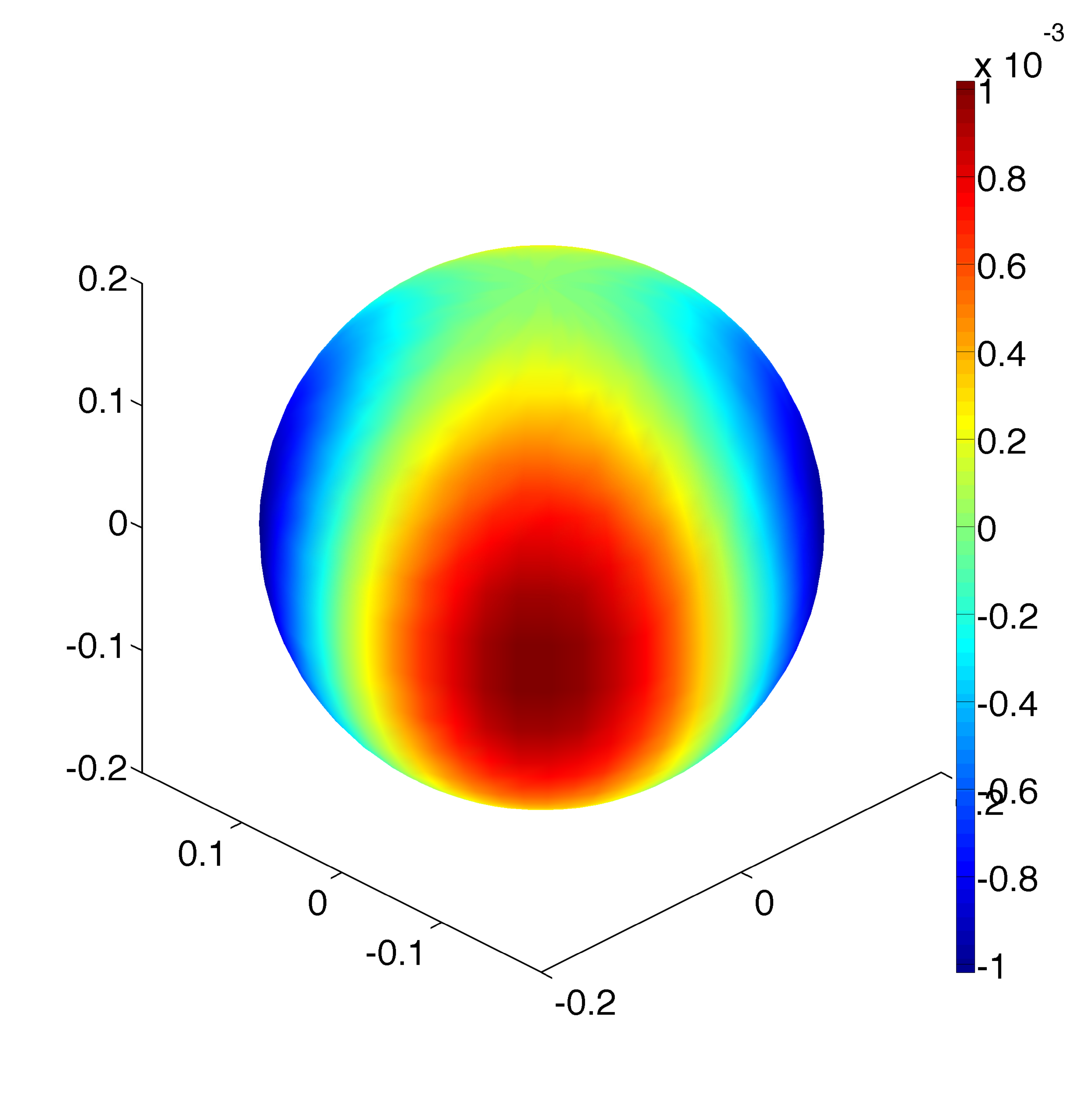} &
\includegraphics[width=0.24\textwidth]{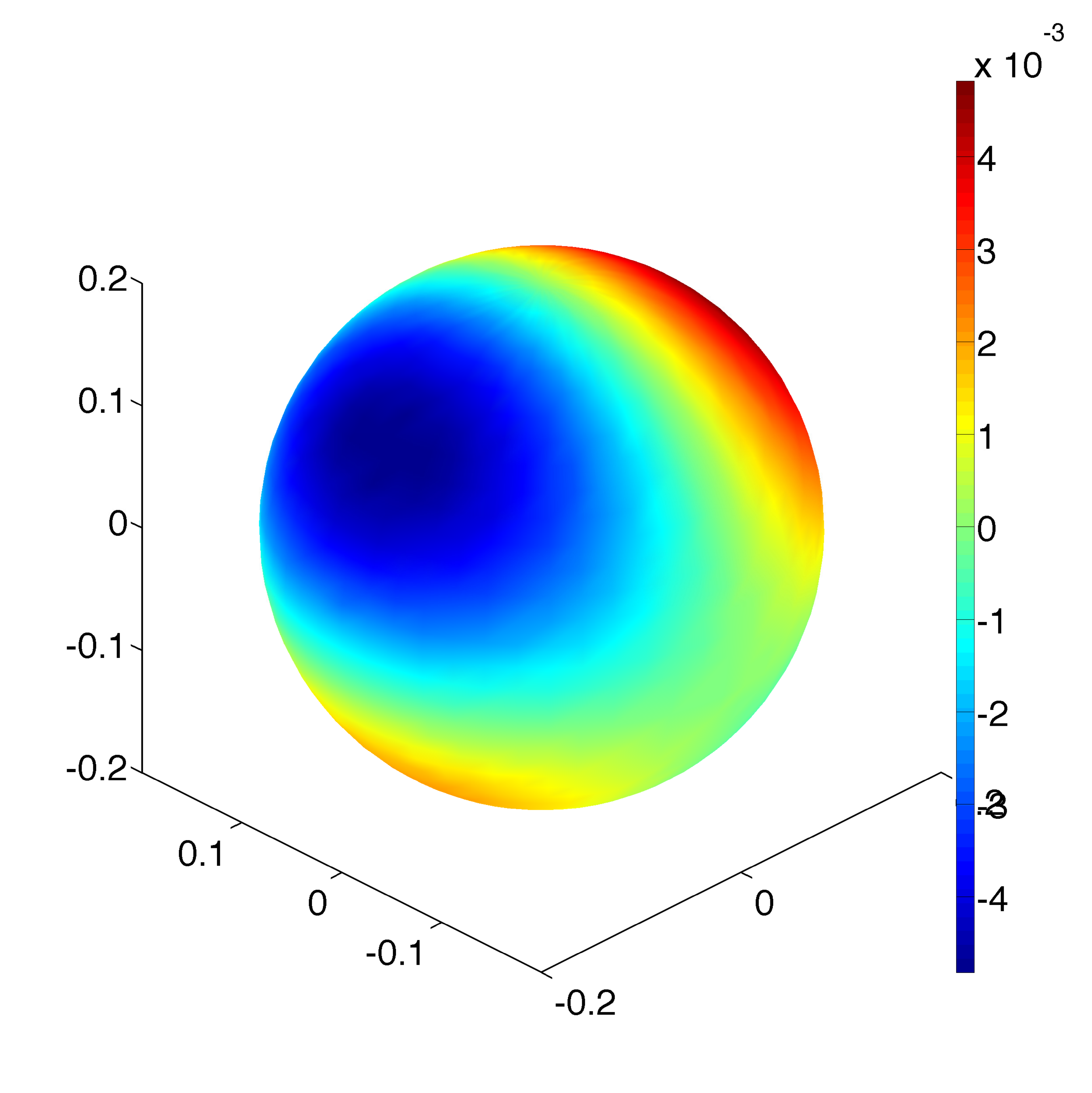}\\
(a) & (b) & (c) & (d)
\end{tabular}
\end{center}
\caption{Numerical solution of scattering of a sphere (radius $a=0.2$) with 42
collocation points. (a) incident wave; (b) $E_{x}$; (c) $E_{y}$, (d)$E_{z}$.}%
\label{fig:VIE-sphere}%
\end{figure}

Fig. \ref{fig:p-conv-sphere} shows the error in terms of the equivalent
collocation points in one direction. In this case, the Mie solution
(\ref{eqn:Mie}) is used as the reference solution with the parameters $a=0.2$,
$k=1$, $m=\sqrt{2}$ (or $\Delta\epsilon=1$). From the fitted lines, it can be
seen that the $\log_{10}(\mathrm{Error)}$ decays linearly with respect to $p$.
Quantitatively, we have
\begin{equation}
\alpha=-0.38,\quad\beta=-1.95.\nonumber
\end{equation}
Hence, we conclude that the exponential convergence is achieved for the errors
against the number $p$. Fig. \ref{fig:VIE-sphere} plots the numerical solution
of the VIE in the sphere with 42 collocation points.

\begin{itemize}
\item \bigskip Case 2: Solution of the VIE in a cube
\end{itemize}

In the second case, we examine the convergence of the solution in a cube. A
regular Lagrange interpolation is used with Gauss nodes in each direction. As
there is no exact solution, the numerical solution with collocation point
$p=7$ in each direction is chosen as the reference solution $\mathbf{E}%
^{\mathrm{ref}}$. The differences between the reference solution
$\mathbf{E}^{\mathrm{ref}}$ and the solutions $\mathbf{E}^{p},p=3,4,5,6$ are
computed.
%\begin{figure}[ptb]
%\begin{center}%
%\begin{tabular}
%[c]{ccc}%
%\includegraphics[width=0.33\textwidth]{images/conv1} &
%\includegraphics[width=0.33\textwidth]{images/conv2} &
%\includegraphics[width=0.33\textwidth]{images/conv3}\\
%(a) & (b) & (c)\\
%&  &
%\end{tabular}
%\end{center}
%\caption{$\log_{10}$ errors of solutions agains order of polynomial basis funtions for $E_{x}$ (a), $E_{y}$ (b) , and $E_{z}$ (c).}%
%\label{fig:p-conv}%
%\end{figure}
\begin{figure}[ptb]
\begin{center}
\includegraphics[width=0.6\textwidth]{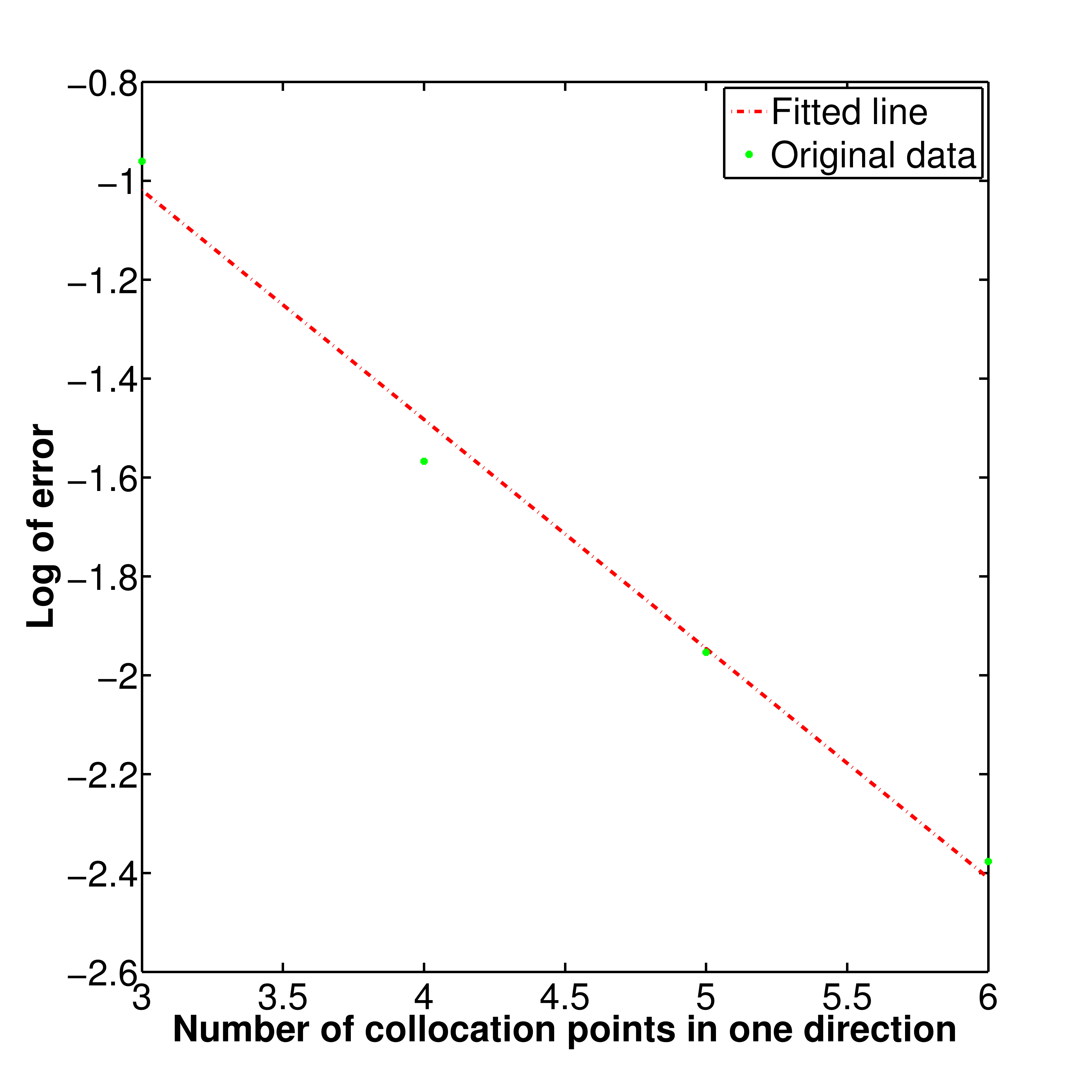}
\end{center}
\caption{$\log_{10}$ errors of solutions $\mathbf{E}^{p}$ in a cube against
number of collocation points in one direction.}%
\label{fig:conv-cube}%
\end{figure}

\begin{figure}[ptb]
\begin{center}%
\begin{tabular}
[c]{cccc}%
\includegraphics[width=0.23\textwidth]{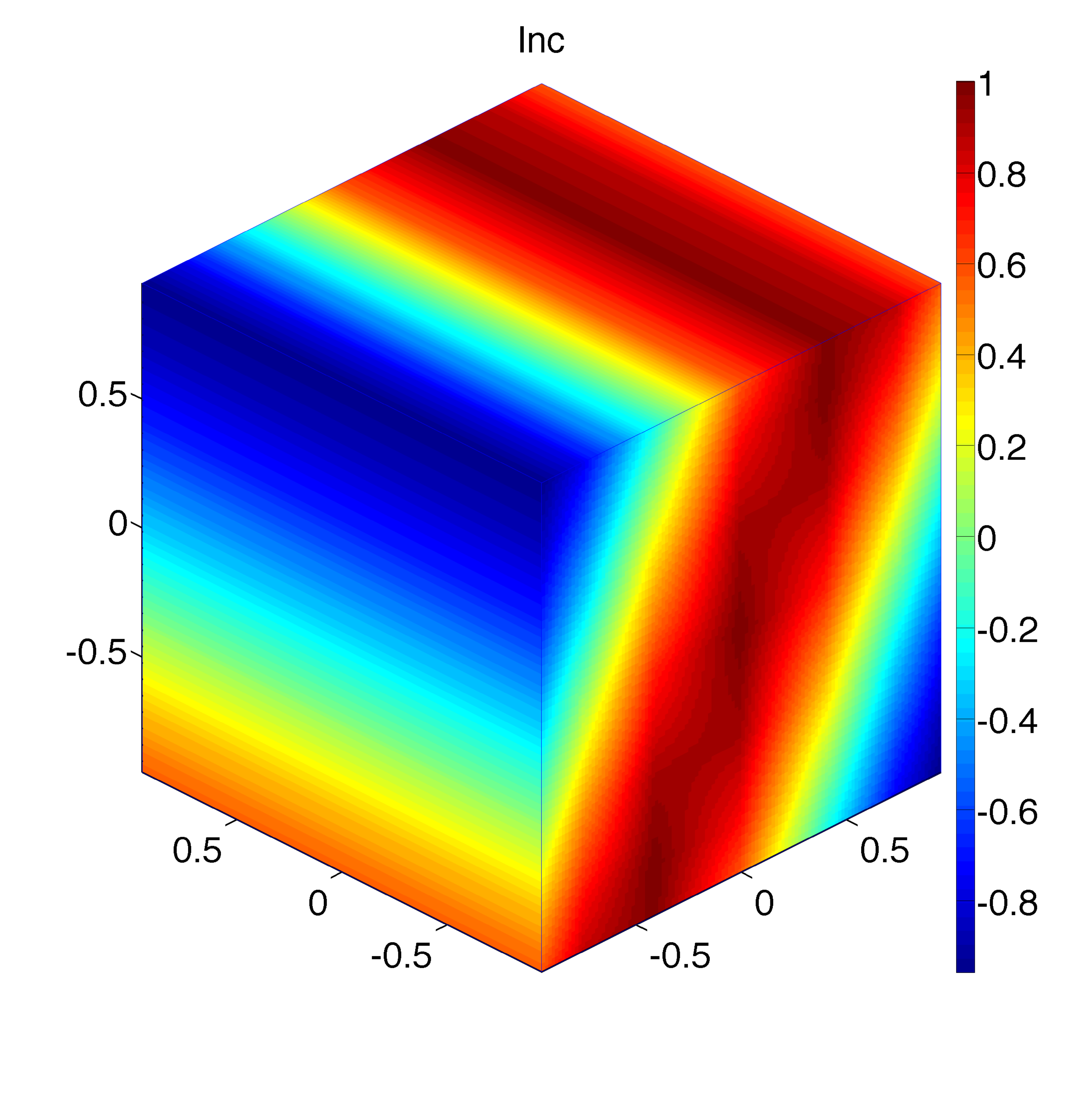} &
\includegraphics[width=0.23\textwidth]{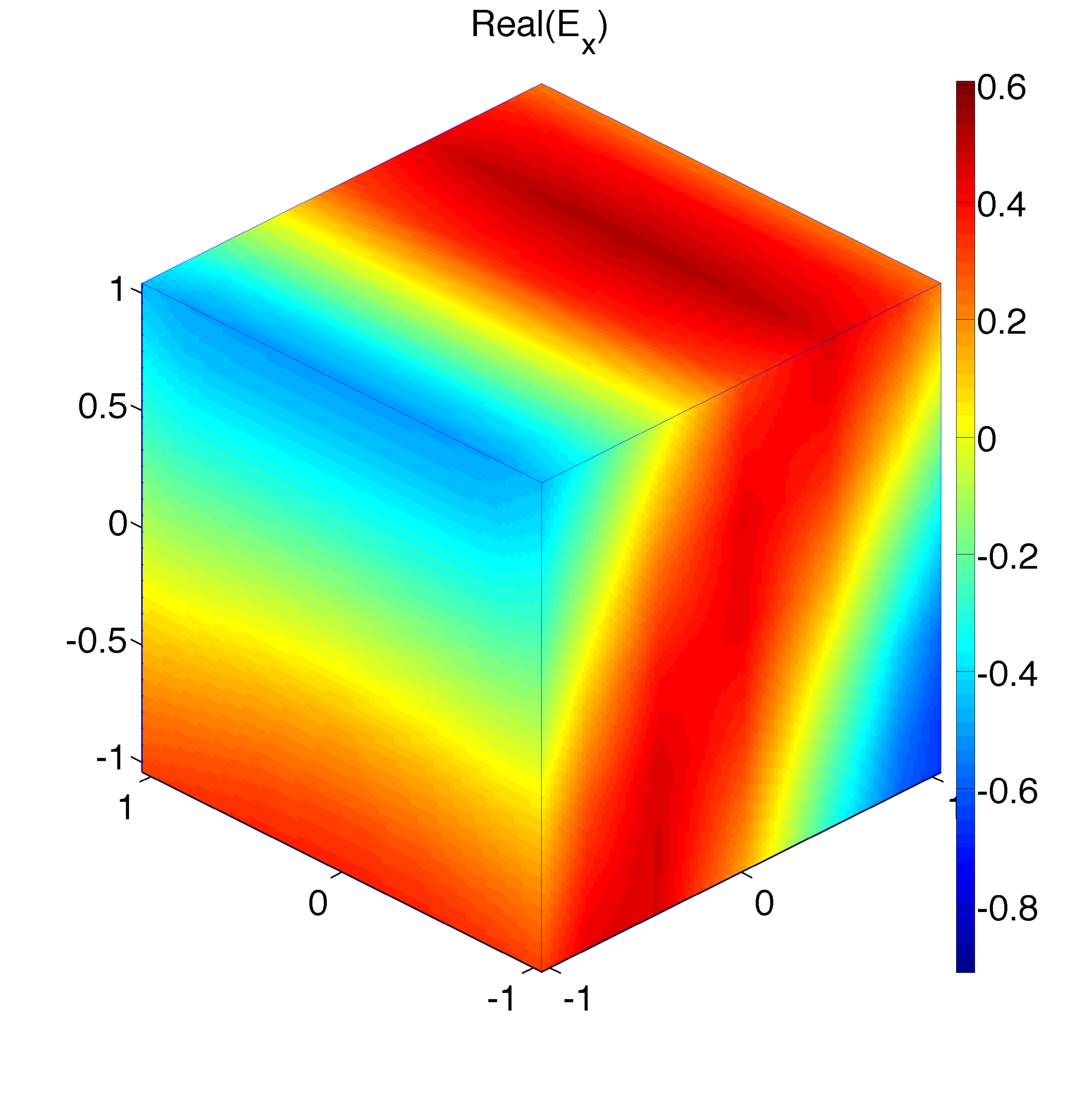} &
\includegraphics[width=0.23\textwidth]{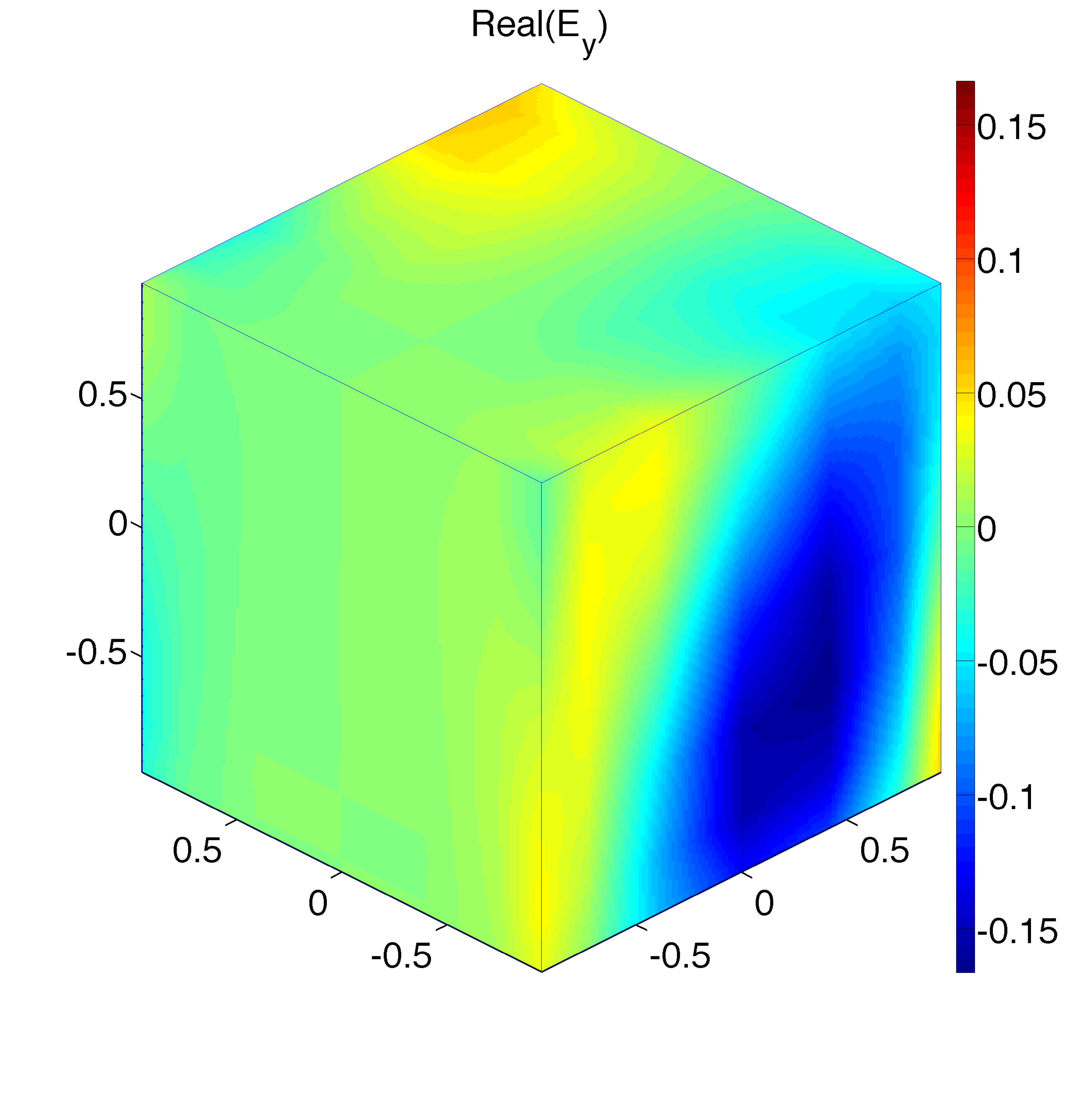} &
\includegraphics[width=0.23\textwidth]{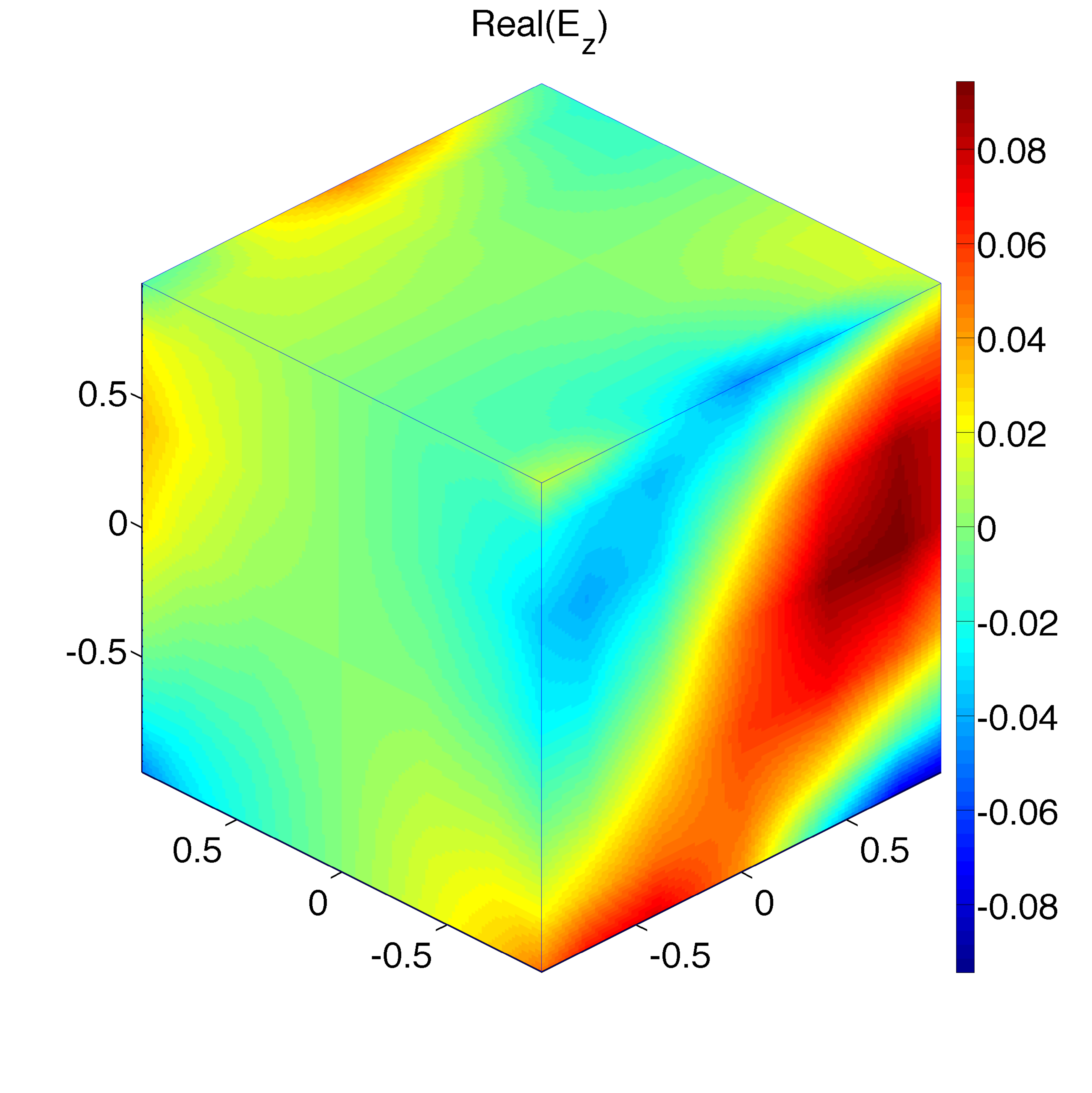}\\
(a) & (b) & (c) & (d)
\end{tabular}
\end{center}
\caption{Incident wave (a) and electric field (b)-(d) for a cubic scatterer
with 27 collocation points.}%
\label{fig:VIE-cube}%
\end{figure}

Fig. \ref{fig:conv-cube} plots the $\log_{10}$ error against $p$, the number
of collocation points along each direction, as well as the fitted line, for
the solution $\mathbf{E}^{p}$. It can be seen that the $\log_{10}%
(\mathrm{Error)}$ decays linearly with respect to $p$. Quantitatively, we
obtain
\begin{equation}
\alpha=-0.464,\quad\beta=-0.092.\nonumber
\end{equation}
Hence, we conclude that exponential convergence is achieved for the errors in
terms of the order $p$. Fig. \ref{fig:VIE-cube} shows the 3-D plot of the
incident wave $\mathbf{E}^{\mathrm{inc}}=\mathbf{i}_{x}e^{ik(-y+0.5z)}$ and
the resulting electric fields over the cube with 27 collocation points.

\subsection{Scattering of multiple scatterers}

Microstructures made of random scatterers, such as the rough surface of solar
cell panels or meta-atoms in meta-materials, can be modeled as an array of
single scatterers of fundamental shapes like cubes and spheres. In this
subsection, we present the capability of our algorithm to handle multiple
scatterers, either in regular or random distributions.

\begin{figure}[ptb]
\begin{center}
%\begin{tabular}
%[c]{ccc}%
\includegraphics[width=0.3\textwidth]{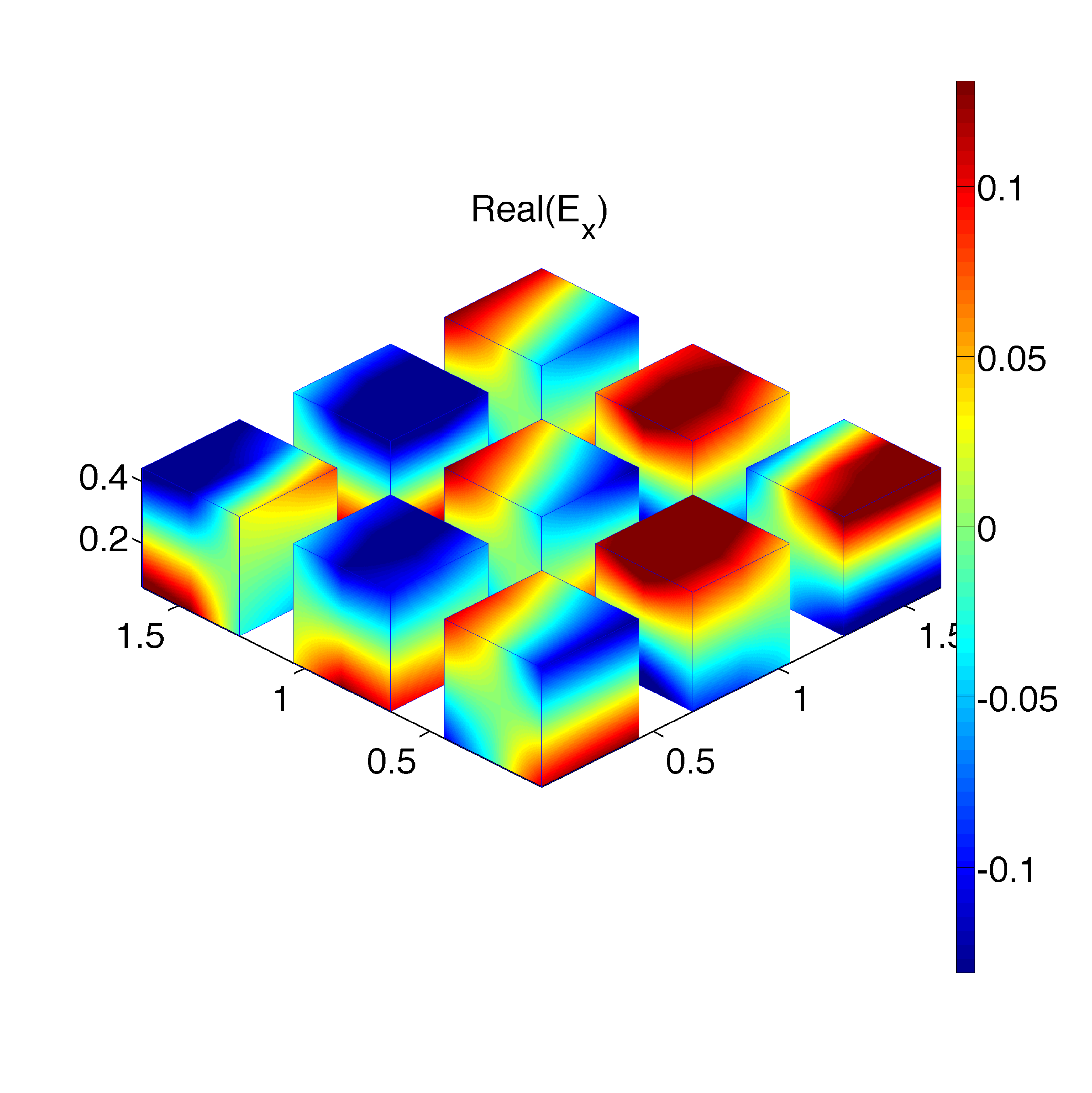}
\includegraphics[width=0.3\textwidth]{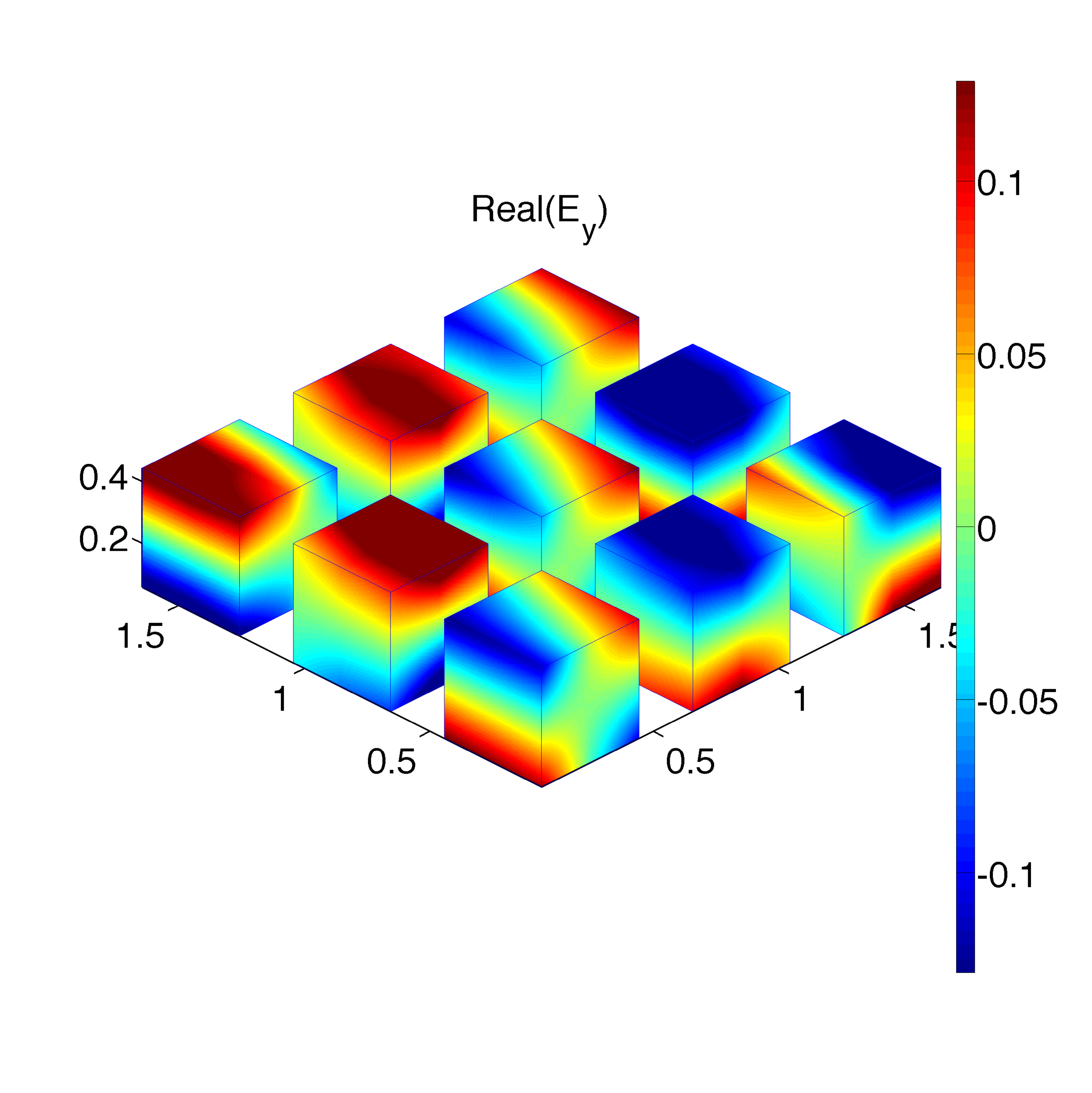}
\includegraphics[width=0.3\textwidth]{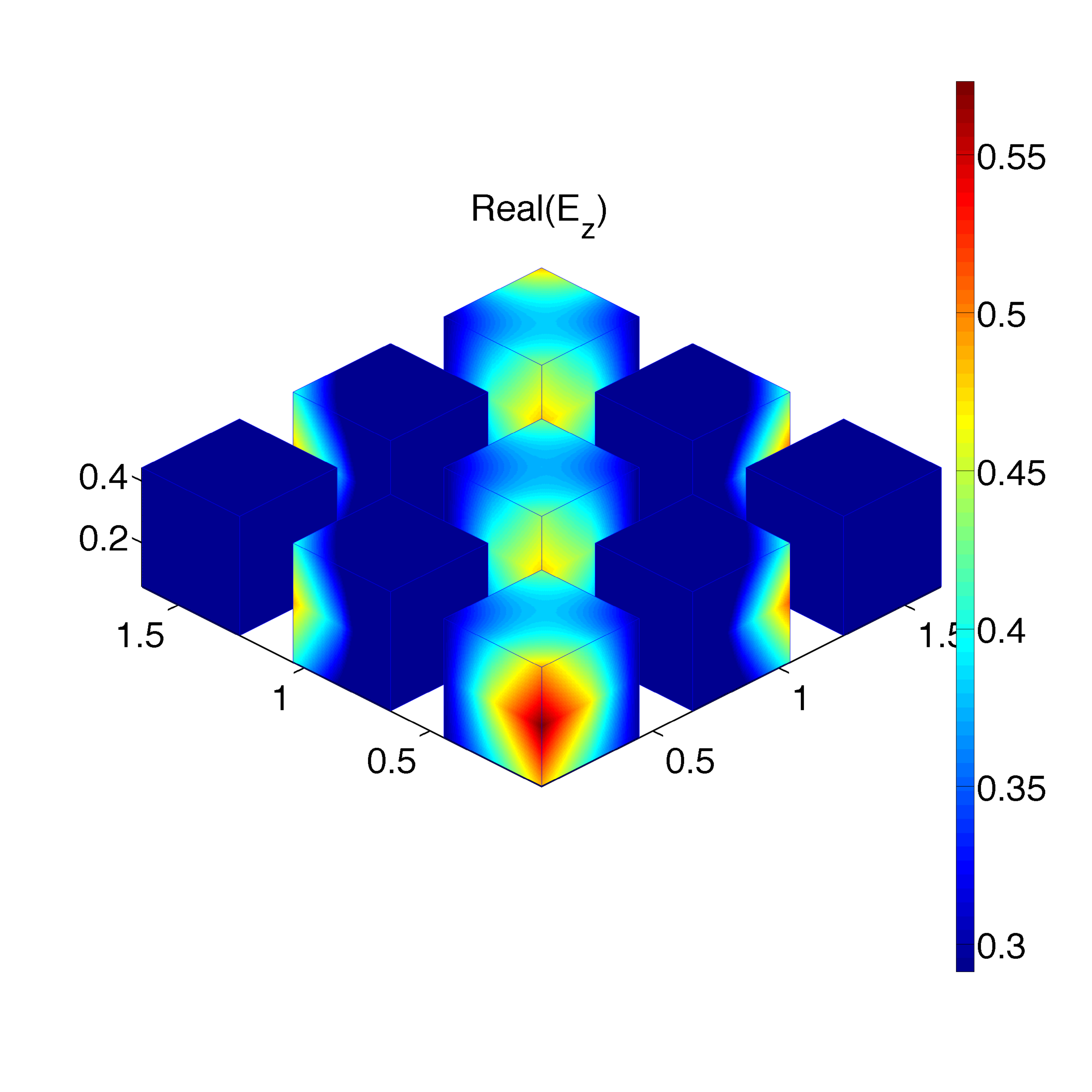}
%\includegraphics[width=0.3\textwidth]{images/ey_cube9}&
%\includegraphics[width=0.99\textwidth]{images/ez_cube9}\\
%(a)& (b) & (c)
%\end{tabular}
\end{center}
\caption{Electric field ($x$-$,y$-, and $z$-components) in a $3\times3$ cube
array}%
\label{fig:cubarrayx}%
\end{figure}

Fig. \ref{fig:cubarrayx} displays the electric field $\mathrm{E}_{x}$,
$\mathrm{E}_{y}$, $\mathrm{E}_{z}$ in the free-space where nine cubic
scatterers are present. In these tests, the incident wave is set as
\begin{equation}
\mathbf{E}_{x}^{\mathrm{inc}}=\mathbf{E}_{y}^{\mathrm{inc}}=0,\quad
\mathbf{E}_{z}^{\mathrm{inc}}=e^{ik(-2x+2y)}. \label{eqn:incube}%
\end{equation}
Each cube has a length of 0.5 and they form a $3\times3$ array align in the
$x$-$y$ plane. The center of the first cube is $(0.25,0.25,0.25)$, and the
remaining cubes are placed 0.1 apart from each other. The parameters are taken
as $\Delta\epsilon=4$ and $\mu=1$. Here, 27 collocation points are used for
each cube.

\begin{figure}[ptb]
\begin{center}
%\begin{tabular}
%[c]{ccc}%
\includegraphics[width=0.3\textwidth]{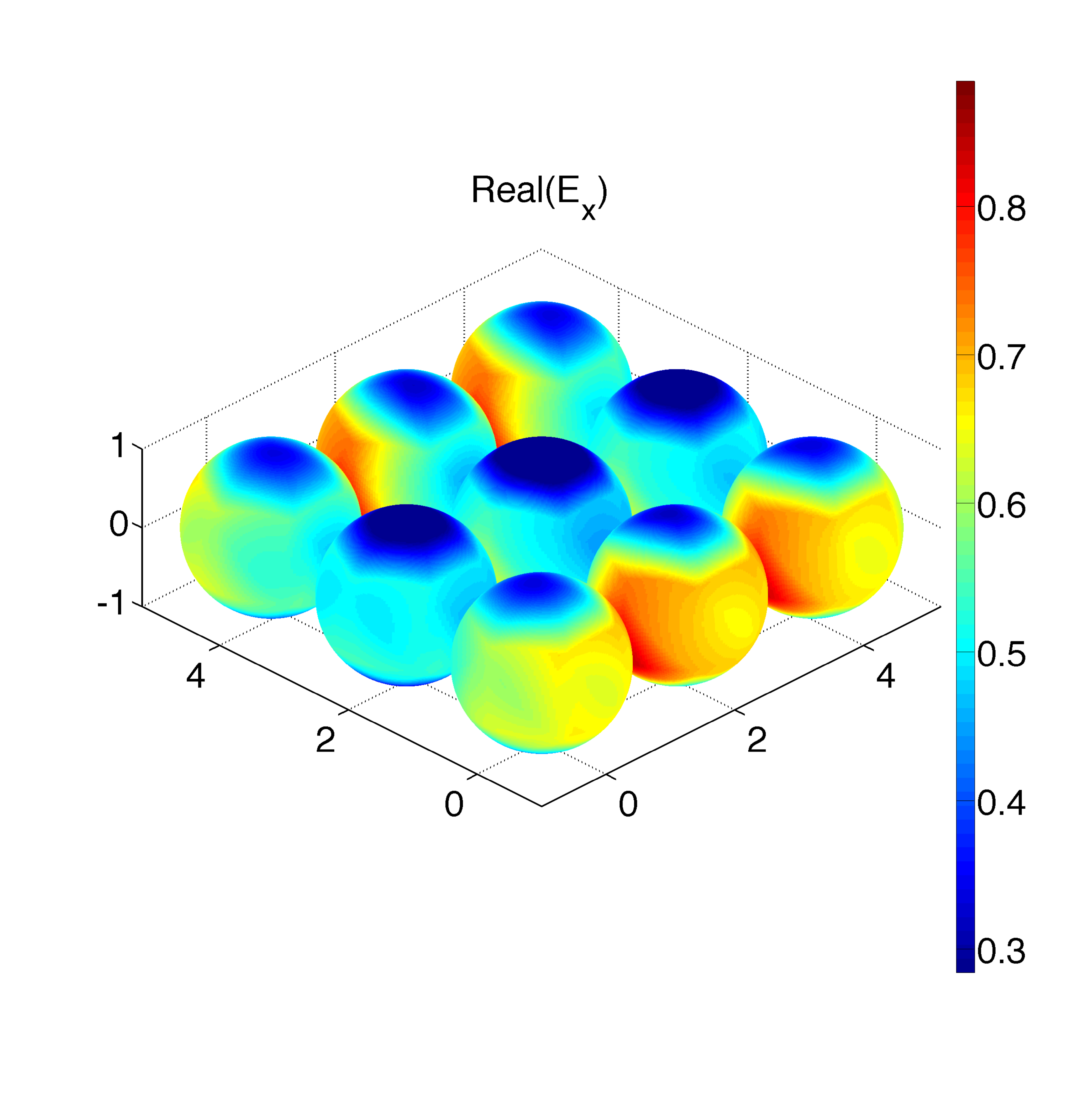}
\includegraphics[width=0.3\textwidth]{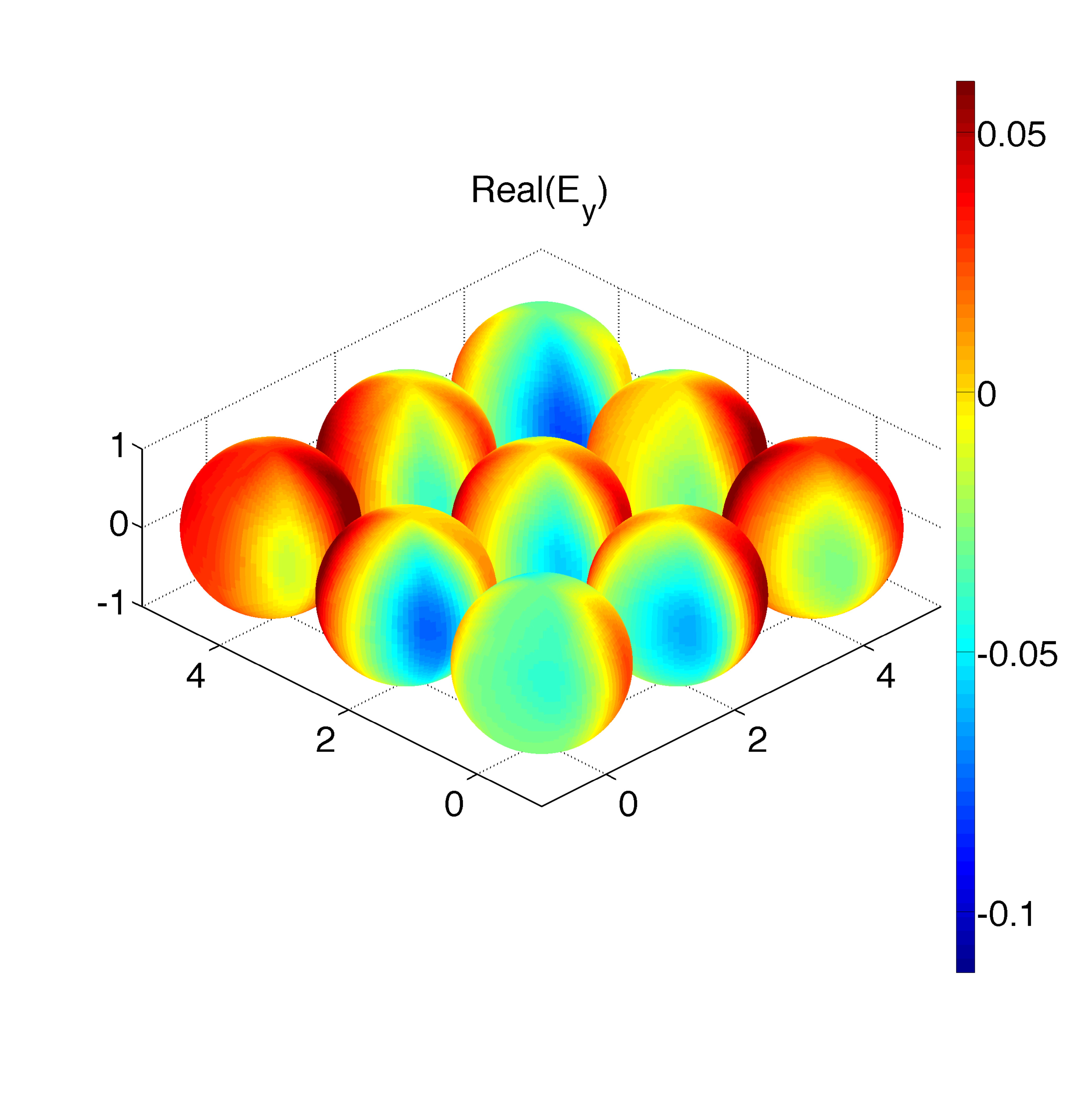}
\includegraphics[width=0.3\textwidth]{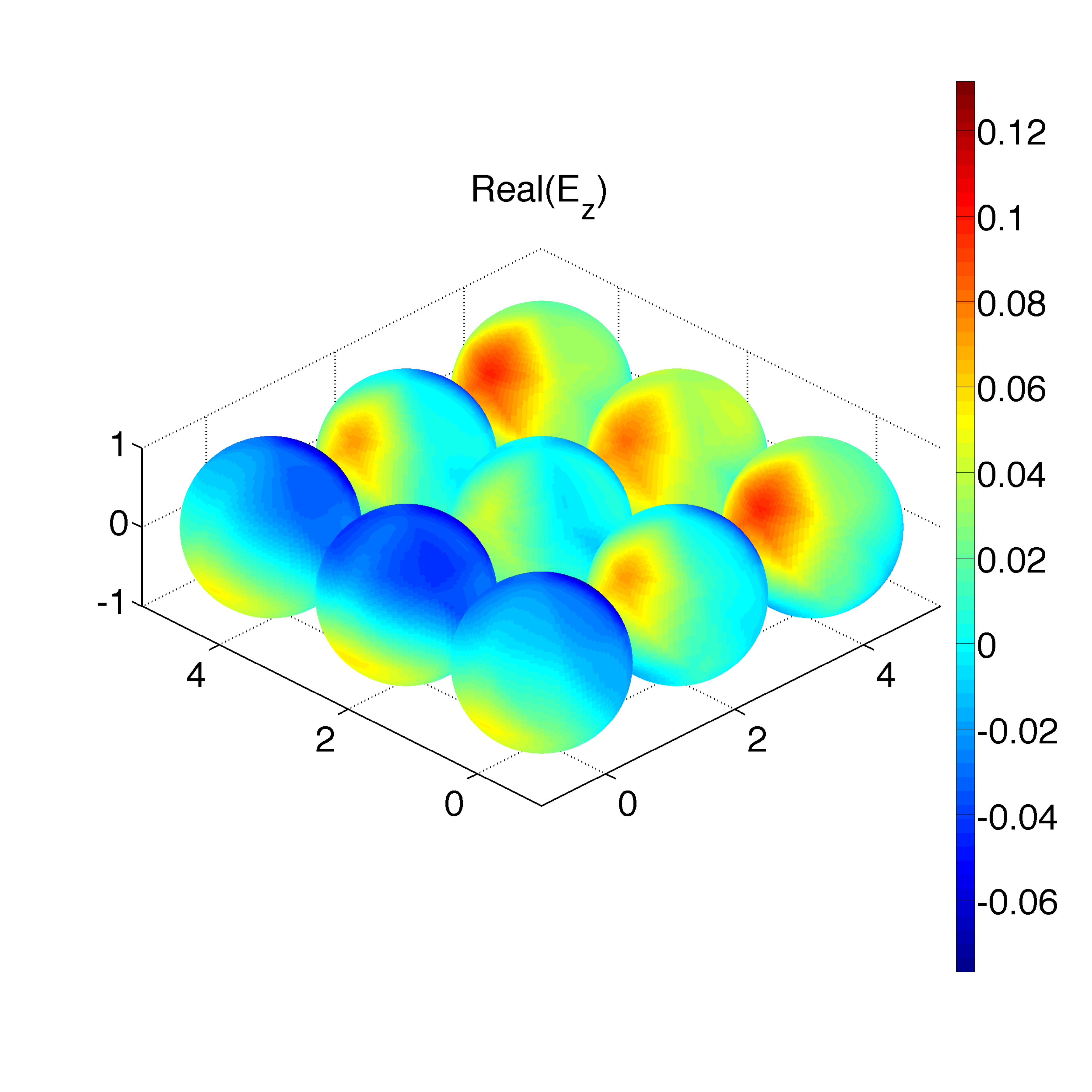}
%\includegraphics[width=0.3\textwidth]{images/ey_cube9}&
%\includegraphics[width=0.99\textwidth]{images/ez_cube9}\\
%(a)& (b) & (c)
%\end{tabular}
\end{center}
\caption{Electric field ($x$-$,y$-, and $z$-components) in a $3\times3$ sphere
array }%
\label{fig:ballarrayx}%
\end{figure}

Fig. \ref{fig:ballarrayx} displays the electric field $\mathrm{E}_{x}$,
$\mathrm{E}_{y}$, $\mathrm{E}_{z}$ in the free-space where nine spherical
scatterers are present. In these tests, the incident wave is set as
\begin{equation}
\mathbf{E}_{x}^{\mathrm{inc}}=e^{ikz},\quad\mathbf{E}_{y}^{\mathrm{inc}%
}=\mathbf{E}_{z}^{\mathrm{inc}}=0. \label{eqn:incball}%
\end{equation}
Each sphere has a radius of 1 and they form a $3\times3$ array align in the
$x$-$y$ plane. The center of the first sphere is $(0,0,0)$ and and the
remaining sphere are placed 0.1 apart from each other. The parameters are
taken as $\Delta\epsilon=1$ and $\mu=1$. Here, 42 collocation points are used
for each sphere.

\begin{figure}[ptb]
\begin{center}
%\begin{tabular}
%[c]{ccc}%
\includegraphics[width=0.45\textwidth]{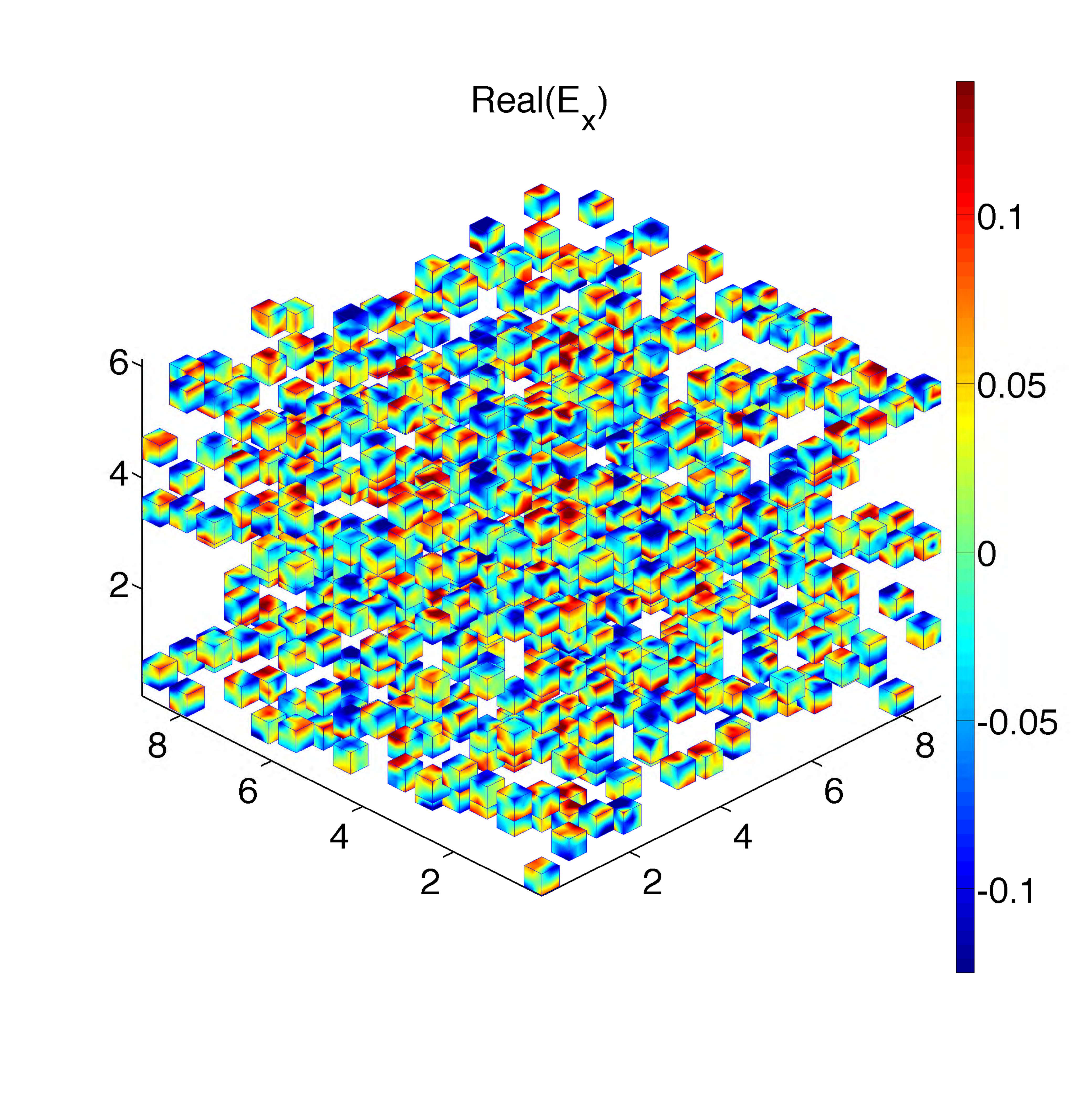}
\includegraphics[width=0.45\textwidth]{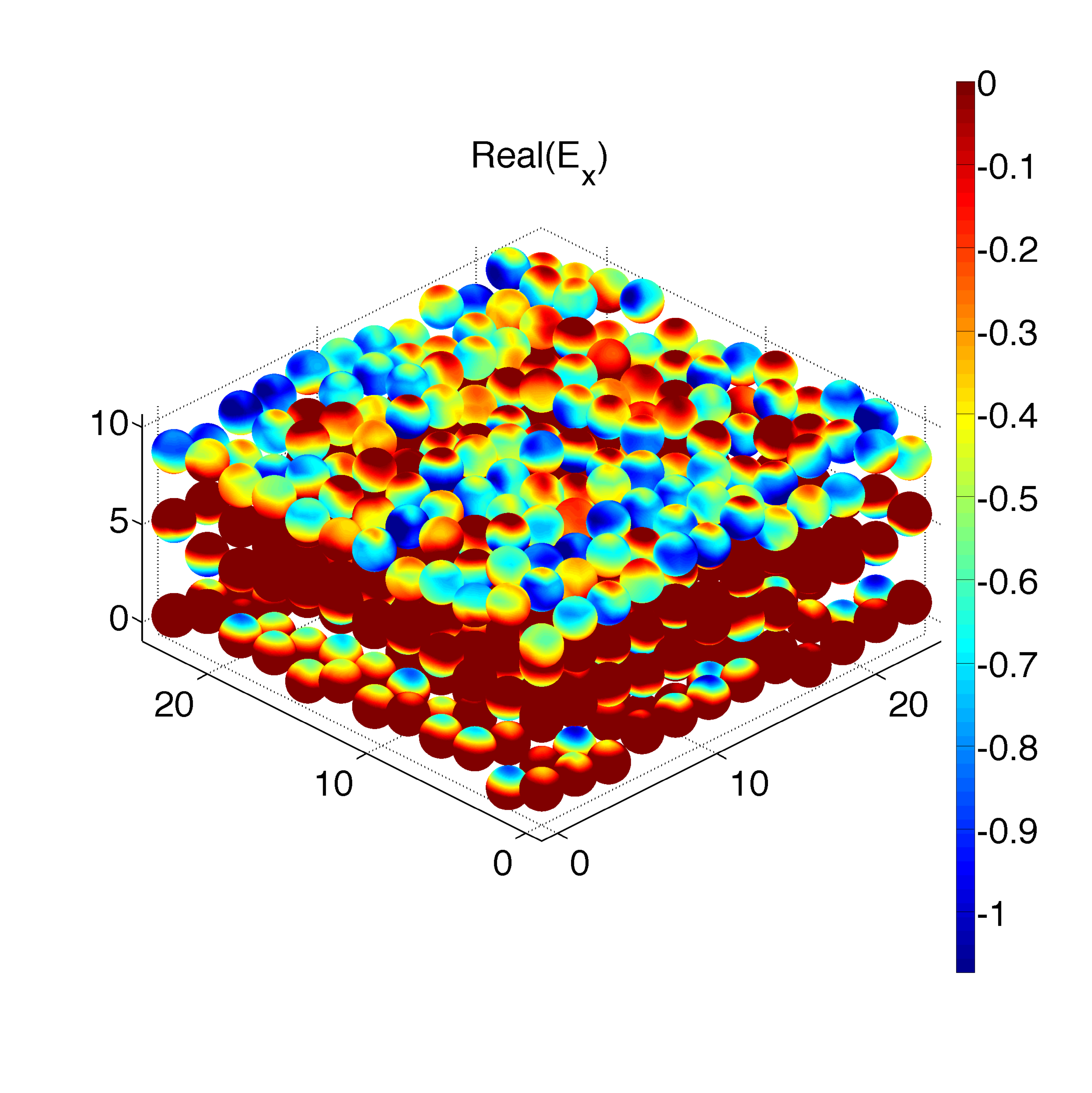}
\par
\includegraphics[width=0.45\textwidth]{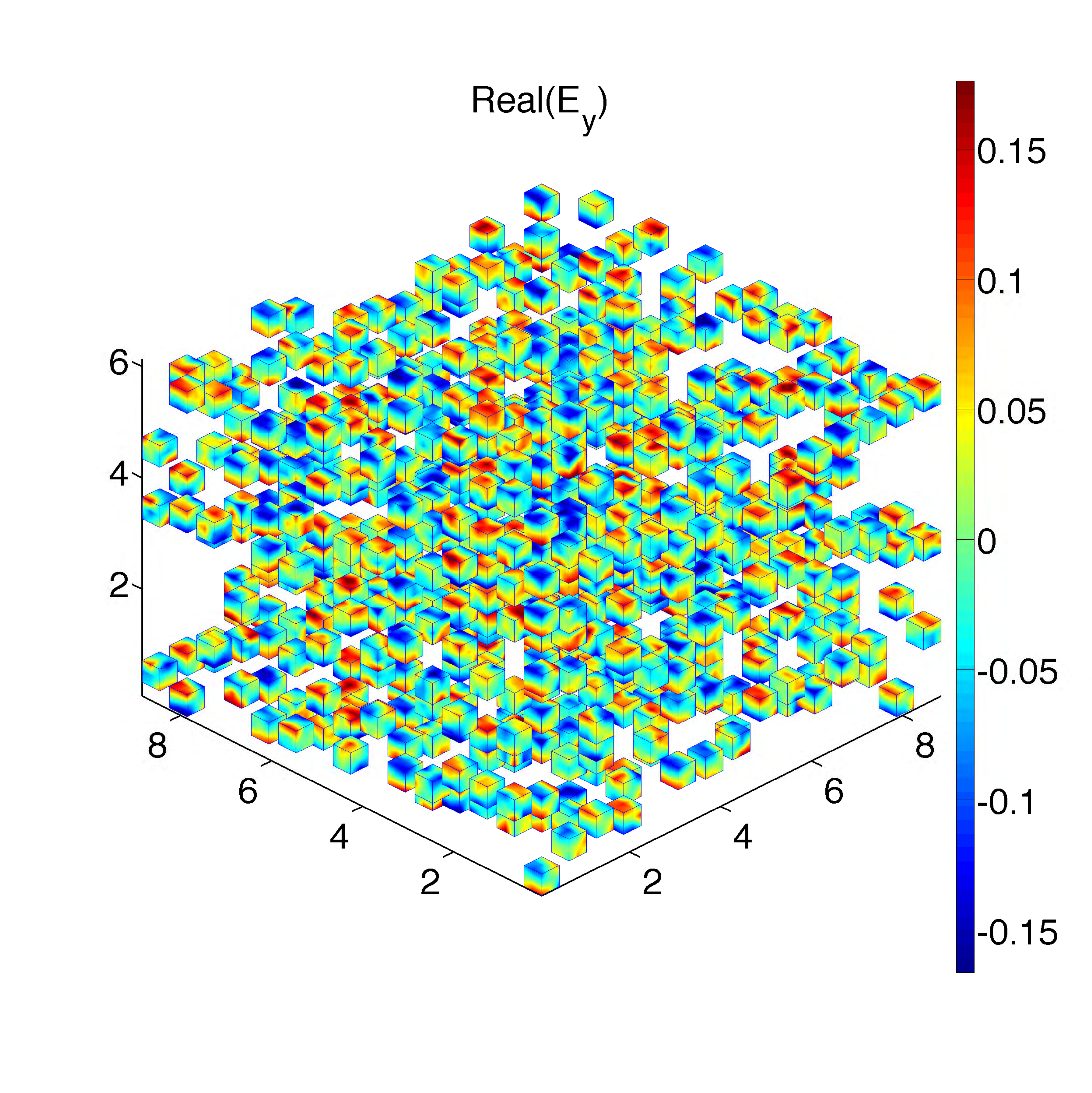}
\includegraphics[width=0.45\textwidth]{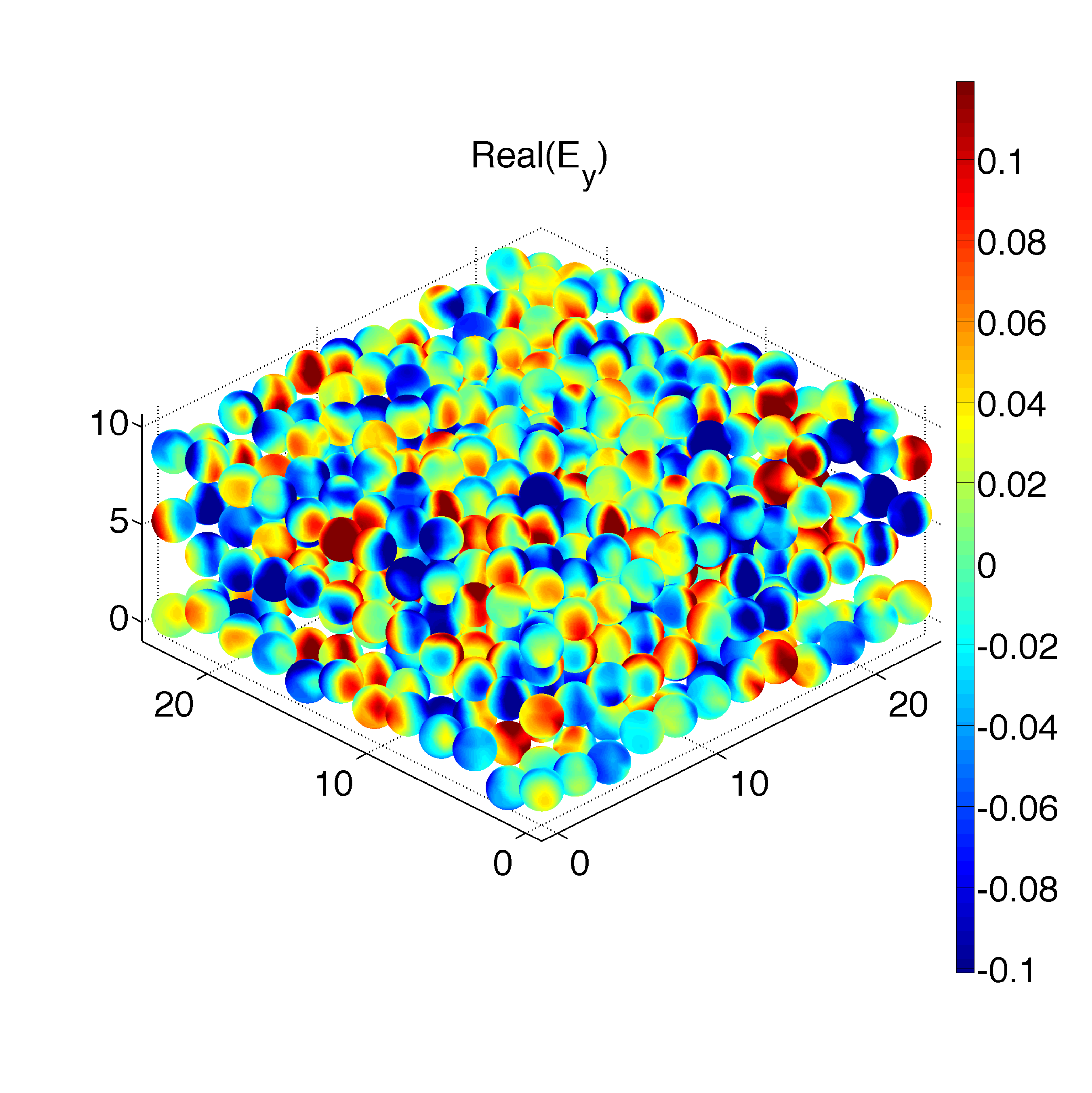}
\par
\includegraphics[width=0.45\textwidth]{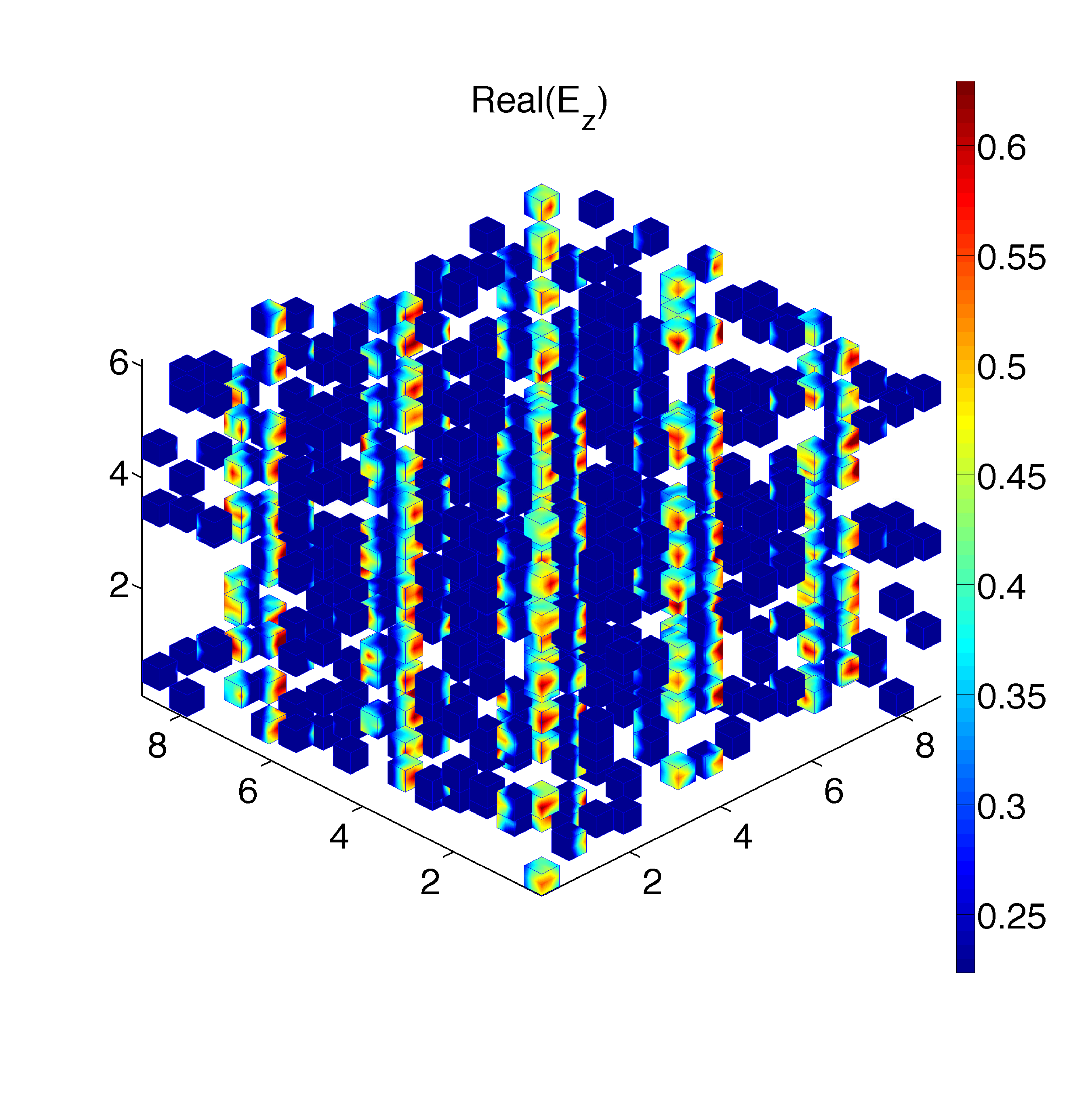}
\includegraphics[width=0.45\textwidth]{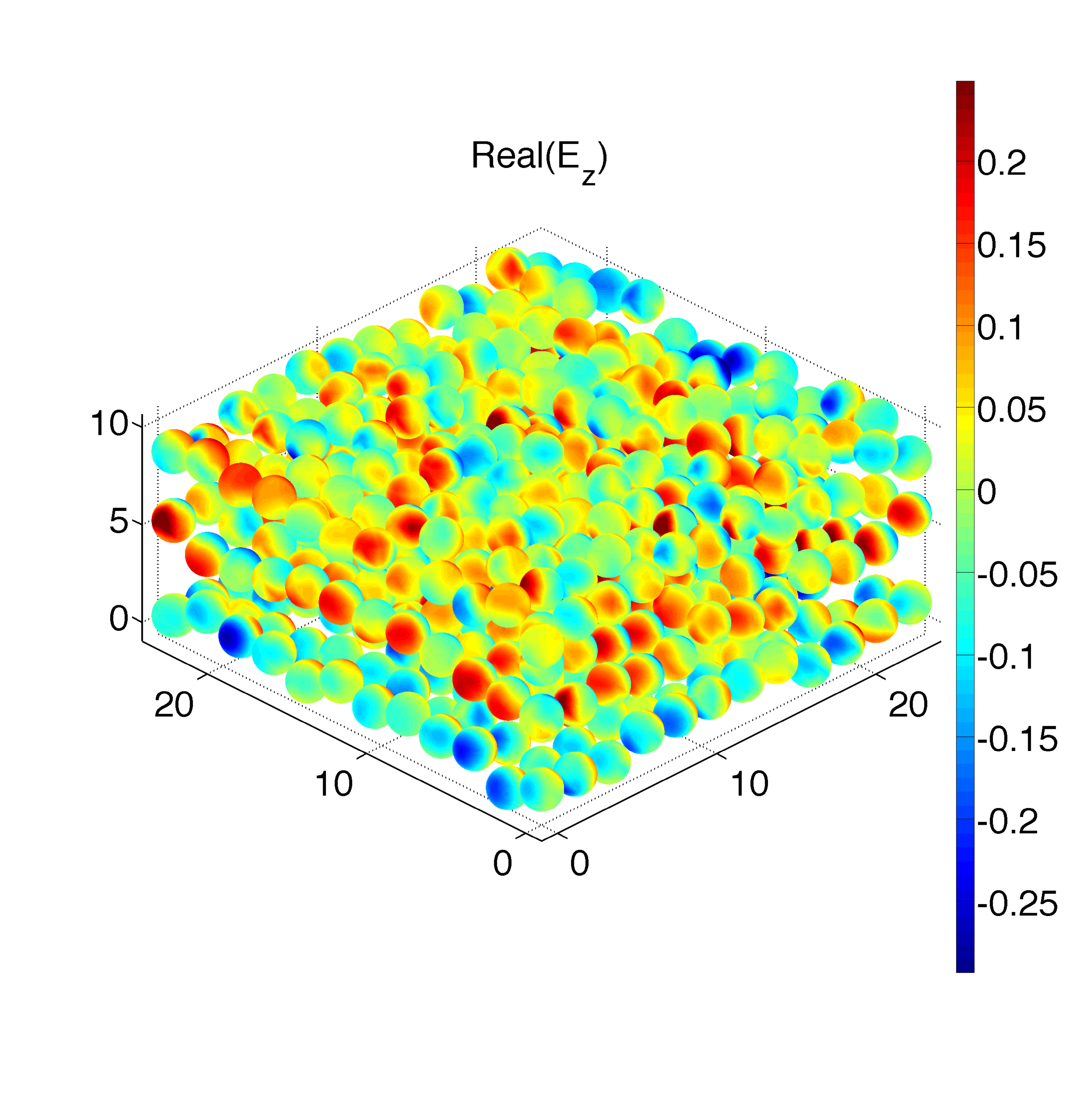}
%\includegraphics[width=0.3\textwidth]{images/ey_cube9}&
%\includegraphics[width=0.99\textwidth]{images/ez_cube9}\\
%(a)& (b) & (c)
%\end{tabular}
\end{center}
\caption{(left) $E_{x},E_{y},E_{z}$ in a random $15\times15\times3$ cube array
with 27 collocation points for each cube of size 0.5; (right) $E_{x}%
,E_{y},E_{z}$ in a random $12\times12\times3$ sphere array with 42 collocation
points for each sphere of radius 1. }%
\label{fig:cube15x}%
\end{figure}

Due to the small number of new quadrature points needed in the Nystr\"{o}m VIE
method, the VIE method can handle hundreds of scatterers of the fundamental
shapes. The left panel of Fig. \ref{fig:cube15x} show the electric field in
675 cubes of size 0.5 under the incident field (\ref{eqn:incube}); these
non-overlapping cubes are arranged in a $15\times15\times3$ random array.
Here, only 27 collocation points are needed for each cube. On the right panel
of Fig. \ref{fig:cube15x}, the electric field in 432 spheres of radius 1 with
the incident field (\ref{eqn:incball}) are shown. Only 42 collocation points
are needed for each sphere and these non-overlapping spheres are arranged in a
$12\times12\times3$ random array. In all the computations, the matrix filling
is the most time consuming part. Therefore, it is computed in parallel using OpenMP.

\section{Conclusion}

\label{sec:conclusion}

In this paper, we have developed an accurate and efficient Nystr\"{o}m method
to simulate the scattering of multiple cubes and spheres, using the volume
integral equation (VIE) for the electric field with the Cauchy principal value
(CPV) of the singular dyadic Green's function. The new formulation allows the
computation of the CPV using a finite size exclusion volume $V_{\delta}$ and
avoiding the usual truncation errors by including missing corrections terms.
As a result, the numerical solution of the VIE is shown to be $\delta
$-independent. In addition, an efficient quadrature formula is employed to
accurately compute the ${\frac{1}{R^{3}}}$ singular integration over the
domain $\Omega\backslash V_{\delta}$. Together, an efficient, accurate, and
$\delta$-independent Nystr\"{o}m collocation method for the VIE is obtained.

In various numerical tests, we demonstrated the convergence of the VIE for
$p$-refinement with increasing order of basis functions (i.e. collocation
points) inside a spherical and a cubic scatterer. Numerical results for the
scattering of multiple scatterers of these two shapes with a small number of
collocation points in each scatterer are also provided.

One of the remaining issues is the treatment of possible field singularities
due to the geometric corners/edges in cubes and cylinders \cite{Bladel} where
specially designed interpolation algorithms may be needed and corresponding
interpolated quadratures will be produced. We will also study the VIE method
for scattering of multiple objects embedded in layered-media for which
treatment of the singularity of layered Green's functions and fast solvers for
the linear system from the Nystr\"{o}m method will be addressed.

\section*{Acknowledgement}

The authors acknowledge the support of the US Army Office of Research (Grant
No. W911NF-14-1-0297).

\bibliographystyle{plain}
\bibliography{vie}

\end{document}